\documentclass[preprint,12pt]{elsarticle}
\usepackage{amssymb}
\usepackage{amsmath}

\usepackage{amsmath}
\usepackage{amsthm}
\usepackage{amssymb}
\usepackage{graphicx}
\usepackage{graphbox}
\usepackage{adjustbox}
\usepackage{xifthen}
\usepackage{algorithm,algpseudocode}
\usepackage{tabularx}
\usepackage{array}
\usepackage{booktabs}
\usepackage{xcolor}
\usepackage{hyperref}
\usepackage{subcaption}

\usepackage[margin=3cm]{geometry}
\makeatletter
\def\ps@pprintTitle{%
  \let\@oddhead\@empty
  \let\@evenhead\@empty
  \def\@oddfoot{}%
  \let\@evenfoot\@oddfoot}
\makeatother

\usepackage{mathtools}
\mathtoolsset{showonlyrefs}

\newcommand{\A}{\mathcal{A}} % arcs/edges in graph
\newcommand{\N}{\mathcal{N}} % nodes in graph
\newcommand{\K}{\mathcal{K}} % set of training data points

\newcommand{\subrefp}[1]{(\subref{#1})}

\newcommand*{\const}[1][]{
\Bigl\lfloor \frac{\kappa_{#1}}{2} \Bigr\rfloor
}

\newcommand{\inputvec}{\nu}
\newcommand{\outputvec}{\hat{x}}
\newcommand{\bipolar}{\{-1,1\}}
\newcommand{\binary}{\{0,1\}}

\newcommand{\bsum}{Y}
\newcommand{\learnrate}{\eta}
\newcommand{\ndrop}{n_{drop}}  % number of dropped neurons in dropout
\newcommand{\minimise}{\textrm{minimise} \ }

\newtheorem{theorem}{Theorem}
\newtheorem{corollary}{Corollary}

\begin{document}

\begin{frontmatter}

\title{Quadratic Unconstrained Binary Optimisation for Training and Regularisation of Binary Neural Networks}

\author[jonas]{Jonas Christoffer Villumsen}
\ead{jonas-christoffer.villumsen@hitachi-eu.com}
\author[sugita]{Yusuke Sugita}
\ead{yusuke.sugita.gt@hitachi.com}

\affiliation[jonas]
            {organization={Research \& Development Group, Hitachi Europe Ltd.},
            city={Kgs.~Lyngby},
            country={Denmark}}

\affiliation[sugita]
            {organization={Research \& Development Group, Hitachi Ltd.},
            city={Kokubunji},
            state={Tokyo},
            country={Japan}}

\begin{abstract}
  Advances in artificial intelligence (AI) and deep learning have raised concerns about its increasing energy consumption, while demand for deploying AI in mobile devices and machines at the edge is growing.
Binary neural networks (BNNs) have recently gained attention as energy and memory efficient models suitable for resource constrained environments; however, training BNNs exactly is computationally challenging because of its discrete characteristics.  Recent work proposing a framework for training BNNs based on quadratic unconstrained binary optimisation (QUBO) and progress in the design of Ising machines for solving QUBO problems suggest a potential path to efficiently optimising discrete neural networks.  In this work, we extend existing QUBO models for training BNNs to accommodate arbitrary network topologies and propose two novel methods for regularisation.
The first method maximises neuron margins biasing the training process toward parameter configurations that yield larger pre-activation magnitudes.   The second method employs a dropout-inspired iterative scheme in which reduced subnetworks are trained and used to adjust linear penalties on network parameters.  We apply the proposed QUBO formulation to a small binary image classification problem and conduct computational experiments on a GPU-based Ising machine.  The numerical results indicate that the proposed regularisation terms modify training behaviour and yield improvements in classification accuracy on data not present in the training set.
\end{abstract}


\end{frontmatter}

\section{Introduction}
Recent advances in deep learning have driven the rapid development of artificial intelligence (AI), enabling widespread deployment in societally critical applications such as autonomous driving, predictive maintenance, and conversational agents.
The emergence of large language models and multimodal generative systems has dramatically expanded the scope of AI applications, yet this progress leads to a substantial rise in energy consumption during both training and inference~\cite{deVries2023}.
Simultaneously, there is a growing demand for integrating AI capabilities into edge environments such as wearable health devices and smart sensors with limited computational and energy resources~\cite{chen2019deep}. 
This dual trend highlights the need for more efficient AI architectures.

To mitigate the computational and energy demands of deep neural networks, quantisation techniques are widely adopted.
Quantisation reduces the precision of network parameters from high-precision floating point representations (e.g., 32-bit) to lower precision (e.g., 8-bit or 4-bit) values, thereby decreasing memory and computational requirements and enabling inference on resource constrained devices.  In the most extreme case, network parameters are reduced to a single bit, resulting in binary (or binarised) neural networks (BNNs) designed for efficient deployment on low-power edge devices~\cite{courbariaux2015, courbariaux2016, hubara2016, rastegari2016}.

The majority of BNN implementations in the literature relies on training through backpropagation algorithms with some adjustments to account for the binarisation of parameters.  In the forward propagation phase the input activations and network parameters are binarised prior to calculating the output activation.  This allows for efficient computation of the pre-activation value using XNOR and popcount operations replacing expensive matrix multiplications~\cite{yuan2023}.
Most of the proposed BNNs use the sign function as a binarisation function; however, it suffers from lack of differentiability and a 0-derivative almost everywhere.
Therefore, when training BNNs using backpropagation and gradient descent, the straight-through estimator is usually applied to approximate gradients through discrete activation functions~\cite{qin2020}.
While this approach enables the practical training of non-differentiable networks, it introduces a gradient mismatch between the forward and backward passes, which can bias the parameter updates and potentially degrade generalisation performance.
Alternatively, model-based approaches, such as mixed-integer and constraint programming formulations~\cite{toro2019}, enable the direct optimisation of parameters without relying on gradient information, but are susceptible to a combinatorial explosion in the number of potential solutions as binary variables are treated endogenously.

Binary neural networks share foundational principles with early energy-based models of neural computation such as the Ising model~\cite{ising1925}, Hopfield networks~\cite{hopfield1982}, and Boltzmann machines~\cite{ackley1985}, which  describe information processing as the minimisation of an energy function characterising the system state.   The Ising model comprises a set of \emph{spins}, each of which can be in one of two states usually represented by the values \(\sigma=1\) or \(\sigma=-1\), pairwise interactions between neighbouring spins governed by a matrix \(J\), and an external field \(h\) acting on individual spins, and where the  Hamiltonian  \(E(\sigma)= - \sum_{i<j }J_{ij} \sigma_i \sigma_j -  \sum_i h_i \sigma_i\) describes the energy of the system.  Hopfield networks describe systems of arbitrarily connected neurons, each with a binary activation (indicating whether to be \emph{firing} or not) that are asynchronously updated based on the weighted sum of inputs exceeding a threshold.  When the strengths of the neuron connections are symmetric, the network can be shown to be isomorphic to an Ising model and thus evolves toward minima of an energy function equivalent to that of an Ising model~\cite{hopfield1982}.  Boltzmann machines generalise Hopfield networks by introducing probabilistic (but still binary) neuron activations to escape from local minima of the energy function~\cite{ackley1985}.
Unlike modern deep neural networks and BNNs trained using backpropagation algorithms, energy-based models do not require that the network topology is acyclic and support arbitrary topologies in principle.

The energy minimisation problem of the Ising model is isomorphic to the problem of solving quadratic unconstrained binary optimisation (QUBO) problems, which provides a general formulation encompassing several classical optimisation problems~\cite{Lucas2014}.
Although QUBO problems are NP-hard~\cite{barahona1982}, a variety of nature-inspired computing systems known as Ising machines have been developed to obtain approximate solutions rapidly and energy-efficiently~\cite{mohseni2022ising}.
Representative implementations include quantum annealers utilising superconducting qubits~\cite{johnson2011quantum}, optical Ising machines based on optical parametric oscillators~\cite{inagaki2016coherent}, and digital annealing systems employing massively parallel architectures such as GPUs or ASICs~\cite{yamaoka201520k,okuyama2019,aramon2019physics,Gotoeaav2372}.
The approach of reformulating constrained optimisation problems as QUBOs amenable to being solved on Ising machines is considered to have the potential to improve the efficiency of combinatorial optimisation, and, hence, Ising machines have been applied to a wide range of practical domains~\cite{Yarkoni_2022}.
In the context of machine learning, applications include the training of Boltzmann machines~\cite{niazi2024training}, clustering~\cite{date2021qubo}, decision tree structure optimisation~\cite{yawata2022qubo}, black-box optimisation~\cite{PhysRevResearch.2.013319,PhysRevResearch.4.023062}, and quantisation of neural networks~\cite{nagel2020}.

QUBO formulations have also been proposed for directly training general neural networks~\cite{abel2022} and BNNs in particular~\cite{sasdelli2021} using quantum annealing.
In particular, Sasdelli and Chin formulate the training problem of layered BNNs as a quadratically constrained binary optimisation problem by ensuring correct neuron activations and minimising the binary training loss, provide a QUBO relaxation, and illustrate the idea on a small network using simulated and quantum annealing~\cite{sasdelli2021}.  Georgiev extends the formulation, presents an iterative algorithm, and performs computational experiments on small instances of image classification and sentiment analysis on a digital Ising machine~\cite{georgiev2023}.  The number of allowed predecessor neurons in the network according to the proposed models is limited to \(2^n-1\) and \(2n-1\), respectively, for some positive integer \(n\).
Since existing methods impose restrictions on the network topology or focus primarily on fitting the training data, it is desired to further investigate methods for incorporating regularisation and generalisation improvements within the QUBO-based BNN framework.

Overfitting arises when network parameters are fitted too closely to the training data and consequently exhibits degraded performance on unseen data.
Regularisation techniques are aimed at improving the model's generalisation performance by preventing overfitting and 
common techniques include L1- and L2-regularisation, early stopping, margin maximisation, and dropout. 
In classification problems, the robustness of a model depends strongly on the geometry of the separating hyperplanes and margin maximisation improves stability by increasing the distance between samples and the decision boundary.  Margin maximisation is a fundamental concept in support vector machines~\cite{cortes1995} and has more recently been applied to  deep neural networks~\cite{sun2015,elsayed2018,jiang2018}.
On the other hand, the dropout technique provides an effective form of regularisation by randomly deactivating nodes during training, which discourages the model from relying excessively on specific features~\cite{hinton2012,srivastava2014}.
Model-based approaches are also susceptible to overfitting~\cite{gambella2021}, and formulations that incorporate regularisation concepts, such as margin maximisation and minimisation of the number of connections, have been proposed for BNNs~\cite{toro2019}.

In this paper we consider the training of BNNs based on a QUBO model amenable to being solved on Ising machines, and introduce two novel approaches for regularisation in QUBO-based BNNs.  The contributions of this work are four-fold.  Firstly, we provide a general QUBO model for training BNNs with arbitrary topology and prove its validity.  Unlike previous work, we do not restrict the number of allowable predecessor nodes.  Secondly, we conduct computational experiments confirming the validity of the model on various network architectures for a four-class image classification problem.  Thirdly, we introduce a novel quadratic penalty expression that maximises neuron margins and improves robustness of the trained networks.  Computational results show significantly increased accuracy of up to 31\% on test data.  Finally, we present a dropout-like iterative algorithm that randomly modifies the network topology leading to modest gains in test accuracy of up to 9\% for some hyper-parameter settings.  All model and algorithmic developments are supported by comprehensive computational experiments on a digital Ising machine.

The remaining part of this paper is laid out as follows.
In Section~\ref{sec:model}, we derive and present the model as a quadratically constrained feasibility problem and its transformation to a QUBO problem.  The experimental set-up is provided in Section~\ref{sec:setup}, where we introduce the image classification problem and details on the computational environment.  In Section~\ref{sec:architectures}, we apply the model developed in Section~\ref{sec:model}, to the classification problem without regularisation using various neural network architectures and present the computational results.  Having established the validity of the model, we introduce in Section~\ref{sec:sum-of-margins}, a new quadratic penalty term that maximises neuron margins and show how it improves accuracy on data not available in the training of the model.  In Section~\ref{sec:dropout}, we present an iterative method that randomly modifies the network topology to further improve test accuracy.  Finally, in Section~\ref{sec:conclusion}, we conclude and discuss further extensions and potential future research directions.

\section{Model formulation}
\label{sec:model}

This study is concerned with BNNs in which all network parameters and activations are either -1 or 1,  and where the  activation function \(f\) is given by,

\begin{align}
f(x) &= \left\{ \begin{array}{rl} 1 & \text{if\ }  x  > 0 \\
                                                  -1 & \text{if\ } x \leq 0 
          \end{array} \right. \label{eq:dummy}
\end{align}

\begin{figure}
  \centering
  
  \begin{subfigure}{0.45\textwidth}
    \centering
    \includegraphics[width=\linewidth]{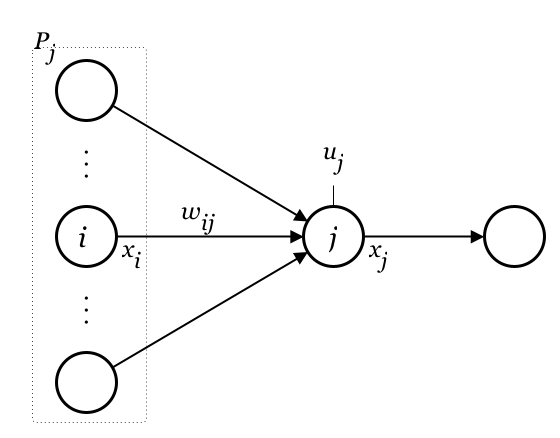}
    \caption{ }
    \label{fig:single-node}
  \end{subfigure}
  \hspace{1cm}
  \begin{subfigure}{0.45\textwidth}
    \centering
    \includegraphics[width=\linewidth]{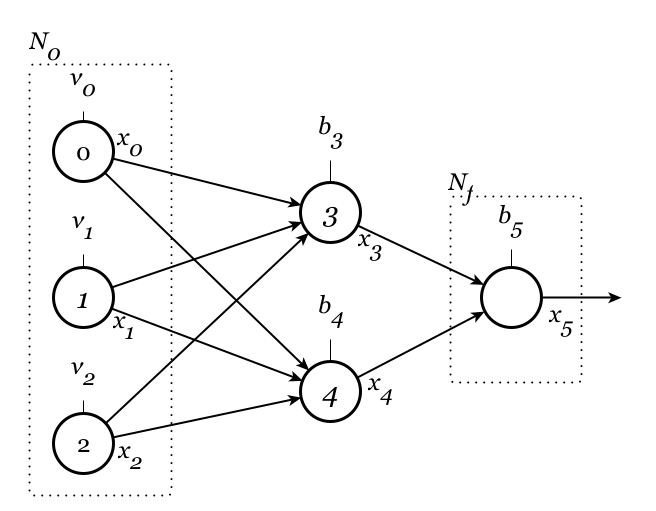}
    \caption{ }
    \label{fig:six-node-network}
  \end{subfigure}

\caption{Schematic diagram of a single node $j$ and its predecessors~\protect\subrefp{fig:single-node} and a small layered network with six nodes~\protect\subrefp{fig:six-node-network}.  
}\label{fig:example-network}
\end{figure}

Let $G=(\N, \A)$ be a directed graph representing the topology of a neural network, such that each node \(i\in\N\) represents a neuron with activation \(x_i \in\bipolar\), and either a fixed input \(\inputvec_i \in \mathbb{R}\) or a bias \(b_i \in \bipolar\), and such that each edge \((i,j)\in \A\) represents a connection between neuron \(i\) and \(j\) with weight \(w_{ij} \in \{-1,1\}\).
Furthermore, let \(P_j=\{i\in\N|(i,j)\in\A\}\) be the set of predecessor nodes of node \(j\) with respect to \(G\) and let \(\N_0\) and \(\N_f\) be the set of input and output nodes, respectively.  Each input node \(i\in\N_0\) has a fixed input \(\inputvec_i \in\mathbb{R}\), while all non-input neurons \(i\in\N \setminus \N_0\) have a variable bias parameter \(b_i \in\{-1,1\}\).
Please refer to Figure~\ref{fig:example-network} where we provide diagrams of single node and a simple network along with the defined notation.
For simplicity of the presentation, we replace the input and bias parameters, for all nodes \(i\in \N\), with a new parameter \(u_i \in \{-1,1\}\), such that \(u_i=f(\inputvec_i)\) if \(i\in \N_0\) and \(u_i=b_i\), otherwise.  For each node \(j \in\N\), the activation of node \(j\) is defined as,

\begin{align}
x_j &= f\left(u_j + \sum_{i\in P_j} w_{ij} x_i \right)  \label{eq:activation}
\end{align}

Please refer to Table~\ref{tbl:notation} where we summarise the symbols introduced here and used throughout the paper.
Note that, we use the terms \emph{bipolar} and \emph{binary} to describe parameters and variables that can take values in \(\{-1,1\}\) and \(\{0,1\}\), respectively.

In the remaining part of this section, we present a QUBO model for training BNNs following a similar approach to the previous studies \cite{sasdelli2021,georgiev2023}.  Contrary to previous work, however, we do not restrict the number of predecessor nodes \(|P_j|\) and we show that our model is valid for any number of predecessors.  Note, that our formulation does not assume a layered architecture and, in fact, connections between any two nodes are allowed in principle, although, in this paper we demonstrate experimental results on standard layered networks only.

In Section~\ref{sec:activation_function}, we encode the activation function into an equality constraint with a binary expansion of the pre-activation whose most significant bit acts as the activation indicator.  In Section~\ref{sec:training_model}, we formulate the training problem by replicating the activation constraints for the data samples and show the transformation to QUBO form.

\begin{table}[htbp]
\caption{\label{tbl:notation}Table of notation.}
\centering
\footnotesize
  \begin{adjustbox}{width=\linewidth}
  \begin{tabularx}{\textwidth}{lcX}
\toprule
Symbol & Domain & Description\\[0pt]
\midrule
\(\N\) &  & Set of nodes in the graph representation of the neural network\\[0pt]
\(\N_0\) &  & Set of input nodes in the graph representation of the neural network\\[0pt]
\(\N_f\) &  & Set of output nodes in the graph representation of the neural network\\[0pt]
\(\A\) &  & Set of connections in the graph representation of the neural network\\[0pt]
\(P_j\) &  & Set of predecessor nodes of node \(j\in \N\)  \(P_j= \{i\in \N \mid (i,j)\in \A \}\)\\[0pt]
\(\K\) &  & Set of training data points\\[0pt]
\midrule
\(x_i\) (\(x^k_i\)) & \bipolar & Bipolar activation of node \(i\in \N\) (for training data point \(k\in \K\))\\[0pt]
\(y_i\) (\(y^k_i\)) & \binary & Binary activation of node \(i\in \N\) (for training data point \(k\in \K\))\\[0pt]
\(\outputvec^k_i\) & \bipolar & Fixed output value of node \(i\in \N_f\) for training data point \(k\in \K\)\\[0pt]
\(\inputvec_i\) (\(\inputvec^k_i\)) & \(\mathbb{R}\) & Input at node \(i\in \N_0\) (for training data point \(k\in \K\))\\[0pt]
\(b_i\) & \bipolar & Bipolar bias of node \(i\in \N \setminus \N_0\)\\[0pt]
\(u_i\) (\(u^k_i\)) & \bipolar & Auxiliary parameter for node \(i\), $u_i=f(\inputvec_i)$ if $i\in \N_0$ and $u_i=b_i$, otherwise (for \(k\in \K\))\\[0pt]
\(d_i\) (\(d^k_i\)) & \binary & Auxiliary parameter for node \(i\), \(d_i= (u_i+1)/2\) (for training data point \(k\in \K\))\\[0pt]
\(w_{ij}\) & \bipolar & Bipolar weight of connection \((i,j)\in \A\)\\[0pt]
\(v_{ij}\) & \binary & Binary weight of connection \((i,j) \in \A\), \(v_{ij}=(w_{ij}+1)/2\)\\[0pt]
\(\psi^k_{ij}\) & \binary & Auxiliary variable used to linearise the activation constraint\\[0pt]
\(\pi_i\) & \(\mathbb{Z}\) & Pre-activation value of node \(i\); sum of bipolar values in Eq.~\eqref{eq:pi}\\[0pt]
\(\rho_i\) & \(\mathbb{Z}_0\) & Pre-activation value of node \(i\); sum of binary values in Eq.~\eqref{eq:rho}\\[0pt]
\(s^{(i)}_l\) ($s^{(k,i)}_l$) & \binary & $l$-th bit of the binary expansion of \(\rho_i\) (for training data point $k\in \K$)\\[0pt]
\(n_j\) & \(\mathbb{Z}_0\) & A non-negative integer that determines the bit-length of the expansion for node \(j\) in Eq.~\eqref{eq:n}\\[0pt]
\(\kappa_j\) & \(\mathbb{Z}_0\) & A non-negative integer that offsets the bit expansion in Eq.~\eqref{eq:kappa}\\[0pt]
\(\chi_j\) (\(\chi^k_j\)) & \(\mathbb{Z}_0\) & Auxiliary non-negative integer variable used in the activation constraints for node \(j\) (training data point \(k\))\\[0pt]
\bottomrule
  \end{tabularx}
  \end{adjustbox}
\end{table}

\subsection{Activation function}
\label{sec:activation_function}

In order to derive the QUBO formulation for the activation constraints~\eqref{eq:activation}, we utilise the following theorem.

\begin{theorem}
\label{thm:1}
Given a non-negative integer $m$ and $x_i\in\{-1,1\}$ for $i=1,\dots,m$, define the non-negative integers
\begin{align}
    n &=  \lfloor \log_2(m) \rfloor, \\ 
    \kappa &= 2^{n+1} - (m+1),
\end{align}
where $\lfloor \cdot \rfloor$ denotes the floor function.

Define bits $s_0,\dots,s_n\in\{0,1\}$ by the $(n{+}1)$-bit expansion
\begin{align}
    \sum_{i=1}^{m} \frac{x_i + 1}{2} + \Big\lfloor \frac{\kappa}{2} \Big\rfloor = \sum_{l=0}^{n} 2^l\, s_l.
    \label{eq:binary-expansion}
\end{align}
Then the sum of the bipolar variables $x_i$ is positive if and only if the most significant bit $s_n$ is $1$, i.e.,
\begin{align}
    \sum_{i=1}^{m} x_i \;\ge\; 1
    \quad\Longleftrightarrow\quad
    s_n=1.
\end{align}
\end{theorem}

\begin{proof}
For notational simplicity, we introduce
\begin{align*}
    \bsum =\sum_{i=1}^{m} \frac{x_i + 1}{2} \,\in\, \{0,\ldots,m\} 
    \quad \text{and} \quad
    \delta = \frac{\kappa}{2} - \Big\lfloor \frac{\kappa}{2} \Big\rfloor \,\in\, \{0,0.5 \}.
\end{align*}

Note that \(0\le \bsum\le m=2^{\,n+1}-1-\kappa\) and thus

\begin{align*}
    0
    \le \bsum + \Big\lfloor \frac{\kappa}{2} \Big\rfloor
    \le 2^{\,n+1} - 1 - \kappa + \Big\lfloor \frac{\kappa}{2} \Big\rfloor
    =   2^{\,n+1} - 1 - \Big\lceil \frac{\kappa}{2} \Big\rceil
    \le 2^{\,n+1}-1,
\end{align*}
so the $(n{+}1)$-bit expansion in Eq.~\eqref{eq:binary-expansion} %defining $s_0, \dots, s_n$ 
exists and is unique.

Since $\sum_{i=1}^m x_i = 2\bsum - m$, we have
\begin{align*}
    \sum_{i=1}^m x_i \ge 1
    &\Longleftrightarrow\; 2\bsum - m - 1 \ge 0\\
   &\Longleftrightarrow\; 2\bsum + \kappa - 2^{n+1} \ge 0 \\
    &\Longleftrightarrow\; 2 \left( \bsum + \Big\lfloor \frac{\kappa}{2} \Big\rfloor + \delta - 2^n \right) \ge 0\\
    &\Longleftrightarrow\; \bsum + \Big\lfloor \frac{\kappa}{2} \Big\rfloor \ge 2^n - \delta.
\end{align*}
Because $\bsum + \lfloor \kappa/2 \rfloor$ is an integer and $ 0 \leq \delta < 1$, the last inequality is equivalent to $\bsum + \lfloor \kappa/2 \rfloor = \sum_{l=0}^{n} 2^l s_l \ge 2^n$.
Furthermore, since $0 \le \sum_{l=0}^{n-1} 2^l s_l \le 2^n - 1$, this is equivalent to $s_n = 1$.
\end{proof}

For each node $j \in \N$,
we denote the pre-activation value by $\pi_j$ and its counterpart composed of the corresponding binary terms by $\rho_j$ as follows,
\begin{align}
    \pi_j  &= u_j + \sum_{i\in P_j} w_{ij} x_i,                            \label{eq:pi}  \\
    \rho_j  &= \frac{u_j+1}{2} + \sum_{i\in P_j} \frac{w_{ij} x_i + 1}{2}.  \label{eq:rho}
\end{align}

\begin{corollary}
\label{cor:1}
For some node $j \in \N$ with the set of predecessor nodes $P_j$,
define non-negative integers
\begin{align}
    n_j &= \lfloor \log_2(|P_j|+1) \rfloor, 
    \label{eq:n}
    \\
    \kappa_j &= 2^{n_j+1} - |P_j| - 2,
    \label{eq:kappa}
\end{align}
and bits $s^{(j)}_0,\dots, s^{(j)}_{n_j}$
such that,
\begin{align}
\label{eq:cor}
\frac{u_j+1}{2} + \sum_{i \in P_j} \frac{w_{ij} x_i + 1}{2} = \sum_{l=0}^{n_j} 2^l s^{(j)}_l - \const[j],
\end{align}
then, due to Theorem~\ref{thm:1},
\begin{align}
u_j + \sum_{i\in P_j} w_{ij} x_i  \geq 1 \quad \text{if and only if} \quad s^{(j)}_n=1.
\end{align}
\end{corollary}

For the following formulation, it
is convenient to replace the bipolar variables with their corresponding binary variables such that for every node \(j\in \N\), \(d_j = (u_j+1)/2\), \(y_j=(x_j+1)/2\), and for every connection  \((i,j)\in\A\), \(v_{ij}=(w_{ij}+1)/2\). We can now re-write \(\rho_j\) in terms of binary variables,
\begin{align}
\rho_j & = d_j + \sum_{i\in P_j} \frac{(2v_{ij} -1)(2y_i-1) + 1}{2}\\
       & = d_j +  \sum_{i\in P_j} (2v_{ij} y_i-v_{ij} -y_i + 1).
\label{eq:rho_bin}
\end{align}

In particular, it has been indicated by Corollary \ref{cor:1} that \(\pi_j \geq 1\) if and only if \(s^{(j)}_{n_j} =1\).
Since 
\(y_j =\frac{x_j+1}{2}\) corresponds to the activation of node \(j\),
the restriction
\(y_j=s^{(j)}_{n_j}\) 
is required for consistency.
Hence, 
from Eqs.~\eqref{eq:cor} and \eqref{eq:rho_bin},
by introducing
\begin{align}
\chi_j =
\begin{dcases}
   0 & \text{if } P_j = \emptyset, \\
   \sum_{l=0}^{n_j - 1} 2^l s_l^{(j)} & \text{otherwise }
\end{dcases}
\label{eq:chi}
\end{align}
and replacing $s^{(j)}_{n_j}$ by $y_j$, we obtain

\begin{align}
  d_j +  \sum_{i\in P_j} (2v_{ij} y_i-v_{ij} -y_i + 1)  & = 2^{n_j} y_j + \chi_j^{} - \const[j]. \label{eq:activation-constr2} 
\end{align}

Then, if constraint~\eqref{eq:activation-constr2} is satisfied, \(y_j\) represents the activation of node \(j\).
Note that constraint~\eqref{eq:activation-constr2} also holds for the
case when \(P_j=\emptyset\); we then have that \(n_j=0\) and \(\kappa_j=0\) and therefore \(d_j=y_j\).  
In particular, when \(j\in \N_0\) is an input node, we have that \(y_j = d_j = \frac{u_j+1}{2} = \frac{f(\inputvec_j)+1}{2}\) is fixed.

\begin{figure}[htbp]
\centering
\includegraphics[width=200px]{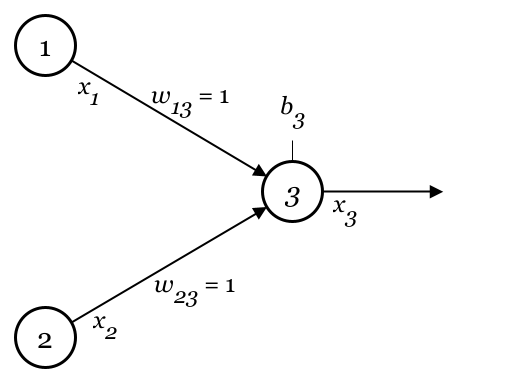}
\caption{\label{fig:example-neuron}Schematic diagram of a neuron \(j=3\) with two predecessor neurons in \(P_j=\{1,2\}\) and bias \(b_3\).}
\end{figure}

Figure~\ref{fig:example-neuron} shows a simple example of the proposed activation encoding using a neuron with two predecessors and all weights equal to 1.
Table~\ref{tbl-activation-example} enumerates all eight combinations of predecessor activations and the neuron bias to demonstrate the exact correspondence between pre-activation, its binary expansion, and the resulting activation.
The pre-activation $\pi_3 = x_1 + x_2 + b_3$ is the sum of three bipolar terms, and its binary counterpart $\rho_3 = y_1 + y_2 + d_3$ counts the number of positive contributions.
The most significant bit in the bit expansion of $\rho_3 = 2s_1^{(3)} + s_0^{(3)}$ uniquely determines whether the bipolar sum is strictly positive.
Indeed, Table~\ref{tbl-activation-example} shows that cases with $\pi_3 \leq 0$ correspond to $\rho_3 \in \{0,1\}$ and $s_1^{(3)} = 0$, whereas cases with $\pi_3 > 0$ correspond to $\rho_3 \in \{2,3\}$ and $s_1^{(3)} = 1$.

\begin{table}[htbp]
\caption{\label{tbl-activation-example}Example of activation and pre-activation variables for a neuron with two predecessor neurons as shown in Figure~\ref{fig:example-neuron} (\(m=3\), \(n=1\), \(\kappa=0\)).  All eight possible combinations of predecessor activations ($x_1$ and $x_2$) and bias ($b_3$) and the corresponding pre-activation ($\pi_3$ and $\rho_3$) and activation values ($x_3$ and $y_3$) are shown.}
\centering
\resizebox{\linewidth}{!}{
\begin{tabular}{ccccccccccc}
\toprule
\(x_1\) & \(x_2\) & \(b_3\) & \(y_1\) & \(y_2\) & \(d_3\) & $\pi_3 = x_1+x_2+b_3$ & \(\rho_3 = y_1+y_2+d_3\) & $s^{(3)}_1 s^{(3)}_0$ & \(y_3 = s^{(3)}_1\) & $x_3$\\[0pt]
\midrule
-1 & -1 & -1 & 0 & 0 & 0 & -3 & 0 & 0 0 & 0 & -1\\
-1 & -1 & ~1 & 0 & 0 & 1 & -1 & 1 & 0 1 & 0 & -1\\
-1 & ~1 & -1 & 0 & 1 & 0 & -1 & 1 & 0 1 & 0 & -1\\
-1 & ~1 & ~1 & 0 & 1 & 1 & ~1 & 2 & 1 0 & 1 & ~1\\
~1 & -1 & -1 & 1 & 0 & 0 & -1 & 1 & 0 1 & 0 & -1\\
~1 & -1 & ~1 & 1 & 0 & 1 & ~1 & 2 & 1 0 & 1 & ~1\\
~1 & ~1 & -1 & 1 & 1 & 0 & ~1 & 2 & 1 0 & 1 & ~1\\
~1 & ~1 & ~1 & 1 & 1 & 1 & ~3 & 3 & 1 1 & 1 & ~1\\
\bottomrule
\end{tabular}
}
\end{table}

\subsection{Training model}
\label{sec:training_model}

In the training model, data-indexed symbols are introduced for a given training set $\mathcal{K}$ (see Table~\ref{tbl:notation}).
Each data point $k\in\mathcal{K}$ consists of an input vector $\inputvec^k$ and a corresponding output vector $\outputvec^k$.
For the data set, the task is to determine weights $w\in\{-1,1\}^{|\mathcal{A}|}$ and biases $b\in\{-1,1\}^{|\mathcal{N}\setminus\mathcal{N}_0|}$.
Hereafter we call an assignment $(w,b)$ \emph{feasible} if, for every $k\in\mathcal{K}$ and every node $j\in\mathcal{N}$, the activation constraints are satisfied.
In addition, for every output node $j\in\mathcal{N}_f$, the activation is required to match the training label, i.e., $x^k_j=\outputvec^k_j$ for all $k\in\mathcal{K}$.
For notational convenience, let $d^k_j$ denote the binarised input/bias parameter for node $j$ and data point $k$; we set $d^k_j=(f(\inputvec^k_j)+1)/2$ if $j\in\mathcal{N}_0$ and set $d^k_j=(b_j+1)/2$ otherwise.

With the above notations, the activation constraints become
\begin{align}
d^k_j + \sum_{i\in P_j} \bigl(2\,v_{ij} y^k_i- v_{ij} - y^k_i + 1\bigr)
= 2^{n_j} y^k_j + \chi^k_j - \const[j],
\quad \forall j\in\mathcal{N},\ k\in\mathcal{K},
\label{eq:training-activation-constraint1}
\end{align}
where, similar to Eq.~\eqref{eq:chi}, for each $j\in\mathcal{N}$ and $k\in\mathcal{K}$ we define $\chi^k_j=0$ if $P_j=\emptyset$ and $\chi^k_j=\sum_{l=0}^{n_j-1} 2^l s^{(k,j)}_l$ otherwise, with $s^{(k,j)}_l\in\{0,1\}$ being auxiliary binary variables.
Here, $y^k_j$ and $v_{ij}$ are binary decision variables corresponding to activations and weights, respectively.
The output consistency is enforced by
\begin{align}
x^k_j = \outputvec^k_j
\quad \Longleftrightarrow \quad
y^k_j = \frac{\outputvec^k_j+1}{2},
\quad \forall j\in\mathcal{N}_f,\ k\in\mathcal{K}.
\label{eq:training-output}
\end{align}

Since constraint~\eqref{eq:training-activation-constraint1} is quadratic,
and due to convenience when formulating the QUBO problem,
we linearise the quadratic terms by introducing for each connection $(i,j)\in \A$ and each training datapoint $k\in \K$ an auxiliary binary variable $\psi^k_{ij}$ and letting

\begin{align}
\psi^k_{ij} = v_{ij}y^k_i, \quad \forall j\in \N, i\in P_j, k \in \K, \label{eq:training-linear-substitution}
\end{align}

We can now define the problem of training a BNN with an underlying topology given by $G=(\N, \A)$ as the non-linear feasibility problem,

\begin{align}
& \textrm{Find vectors} && d \in \binary^{|\N|}, v \in \binary^{|\A|}, y \in \binary^{|\N||\K|}, 
                          \chi \in \mathbb{Z}_0^{|\N||\K|}, \psi \in \binary^{|\A||\K|} \\
&\textrm{such that} &&     d^k_j + \sum_{i\in P_j} \bigl(2\,\psi^k_{ij} - v_{ij} - y^k_i + 1\bigr) = 2^{n_j} y^k_j + \chi^k_j - \const[j],
                         \quad \forall j\in\mathcal{N},\ k\in\mathcal{K}, \label{eq:training-activation-constraint2}\\
&\textrm{and}  &&        \textrm{constraint } \eqref{eq:training-linear-substitution} \\
&\textrm{are satisfied}
\end{align}

To convert the feasibility problem~\eqref{eq:training-linear-substitution}-\eqref{eq:training-activation-constraint2} to QUBO form, we note that a feasible solution to~\eqref{eq:training-activation-constraint2} can be achieved by minimising the quadratic penalty term,

\begin{align}
 H_1  &= \sum_{k\in\K} \sum_{j\in\N}  \left( d^k_j  + \sum_{i\in P_j}  (2 \psi^k_{ij} - v_{ij}-y^k_i+1) - 2^{n_j} y_j^k - \chi^k_j + \const[j] \right)^2
\end{align}

and that constraints~\eqref{eq:training-linear-substitution} can be enforced by the standard penalty \cite{boros2002} given as

\begin{align}
v_{ij} y^k_i -2 (v_{ij} \psi^k_{ij} + y^k_i \psi^k_{ij}) + 3 \psi^k_{ij}, \quad \forall j\in \N, i\in P_j, k\in \K
\label{eq:penalty}
\end{align}

\begin{table}[htbp]
\caption{\label{tbl:binary-substitution-logic}The function \(xy-2(x\psi+y\psi)+3\psi\) is 0 exactly, when \(\psi=xy\) is True and strictly positive otherwise.}
\centering
\begin{tabular}{cccccc}
\toprule
\(x\) & \(y\) & \(\psi\) & \(xy\) & \(\psi=xy\) & \(xy-2(x\psi+y\psi)+3\psi\)\\[0pt]
\midrule
0 & 0 & 0 & 0 & True & 0\\[0pt]
0 & 0 & 1 & 0 & False & 3\\[0pt]
0 & 1 & 0 & 0 & True & 0\\[0pt]
0 & 1 & 1 & 0 & False & 1\\[0pt]
1 & 0 & 0 & 0 & True & 0\\[0pt]
1 & 0 & 1 & 0 & False & 1\\[0pt]
1 & 1 & 0 & 1 & False & 1\\[0pt]
1 & 1 & 1 & 1 & True & 0\\[0pt]
\bottomrule
\end{tabular}
\end{table}

It can be confirmed that constraints~\eqref{eq:training-linear-substitution} hold when the penalty in Eq.~\eqref{eq:penalty} attains its minimum value of 0 as shown in Table~\ref{tbl:binary-substitution-logic}.
We can therefore reformulate the feasibility problem~\eqref{eq:training-linear-substitution}-\eqref{eq:training-activation-constraint2} to that of finding binary vectors $d$, $v$, $y$, $\chi$, and $\psi$ that minimise the following Hamiltonian,

\begin{align}
H_{train}= H_1 + \alpha H_2 \label{eq:training}
\end{align}

where $H_1$ enforces the correct activation of neurons and

\begin{align*}
 H_2 = \sum_{k\in\K} \sum_{j\in\N} \sum_{i\in P_j} v_{ij}y^k_i - 2(v_{ij}\psi^k_{ij} + y^k_i \psi^k_{ij})+3\psi^k_{ij}
\end{align*}

ensures that \(\psi^k_{ij}= v_{ij}y^k_i\)
and $\alpha$ is an appropriate scaling factor.
A feasible assignment of weights and biases are characterised by a solution in which $H_{train} = H_1=H_2=0$ corresponding to satisfying all activation constraints~\eqref{eq:training-activation-constraint2} and the auxiliary constraints~\eqref{eq:training-linear-substitution}.
Complying with constraints in Eq.~\eqref{eq:training-output} is achieved by simply fixing the output variables.

We remark on several practical aspects of the formulation.
First, instead of fixing the output values
as in Eq.~\eqref{eq:training-output}, one could add a loss function to \(H_{train}\) minimising the difference between the output activations and the training output.
However, this would increase the number of variables and terms in the Hamiltonian potentially making it more difficult to find near-optimal solutions.
Second, note that the numbers of variables and constraints scale with the size of the dataset since
for each training data point \(k \in \K\), we introduce variables such as \(y_j^k\), \(\psi_{ij}^k\), and  $\chi^k_j$ for each node \(j\in\N\) and each connection \((i, j) \in \A\).
Therefore, the numbers of variables and constraints grow as \(\mathcal{O}(|\K|\cdot|\N|)\) or \(\mathcal{O}(|\K|\cdot|\A|)\).

\section{Experimental set-up}
\label{sec:setup}

\begin{figure}
  \centering
  \begin{subfigure}{0.17\textwidth}
    \centering
    \includegraphics[width=\linewidth]{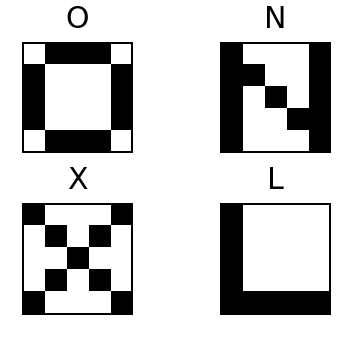}
    \caption{Training images.}
    \label{fig:train_images}
  \end{subfigure}
  \hspace{2cm}
  \begin{subfigure}{0.55\textwidth}
    \centering
    \includegraphics[width=\linewidth]{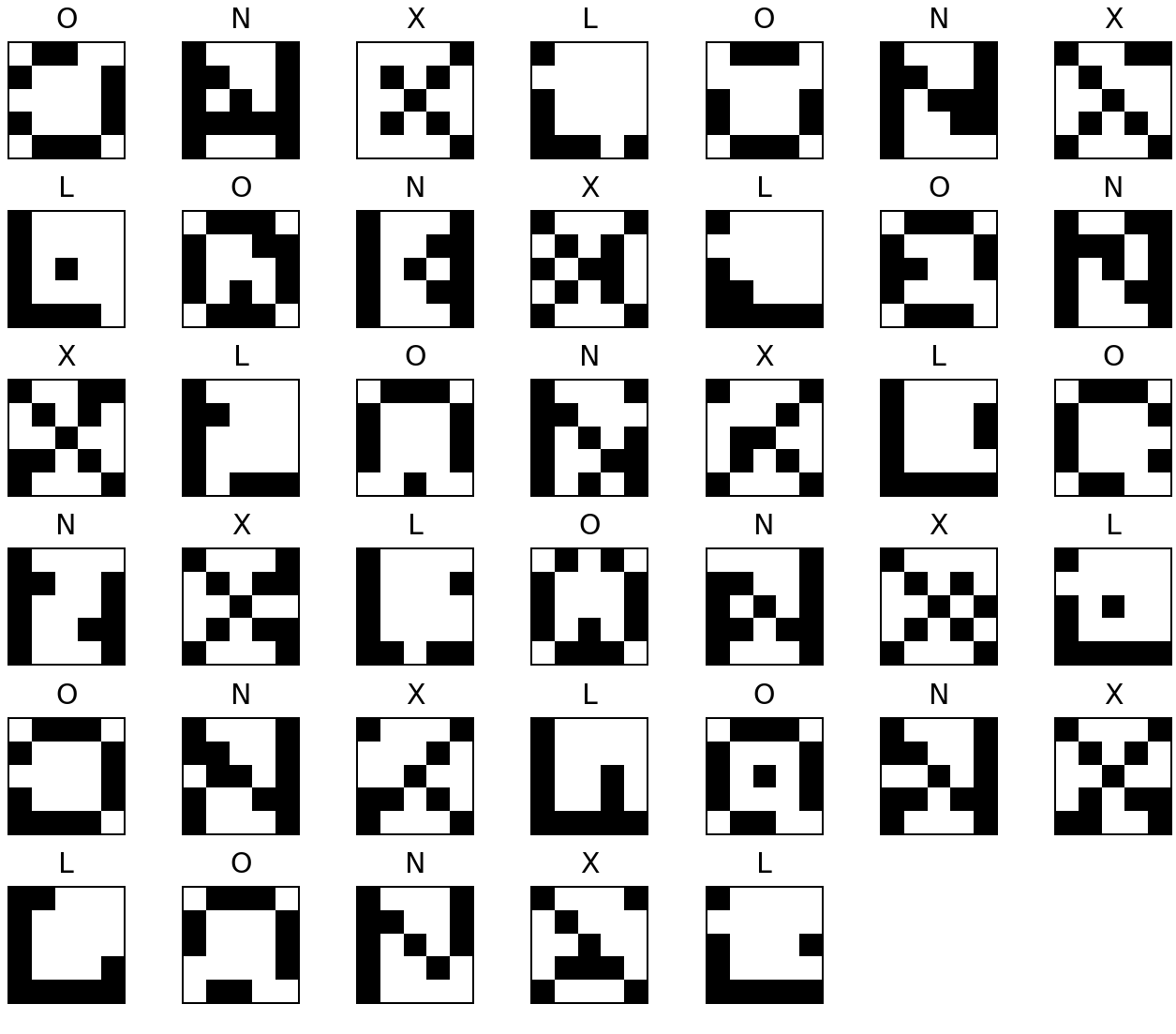}
    \caption{Test images.}
    \label{fig:test_images}
  \end{subfigure}

\caption{Data set comprising 44 monochrome images each labelled with their respective class (O, N, L, or X).  Four images~\subrefp{fig:train_images} are used for training, while 40 images~\subrefp{fig:test_images} are used for testing.}
\label{fig:dataset}
\end{figure}

We apply the training model developed in Section \ref{sec:model} to a four-class image classification problem using one-shot training, where network parameters are learned from a single image from each class.  The data set consists of 44 binary images each comprising \(5\times5\) pixels resembling the letters \emph{O}, \emph{N}, \emph{L}, and \emph{X}.   Each image is classified accordingly into one of four classes O, N, L, and X.  The data set is shown in Figure~\ref{fig:dataset} and contains eleven images from each class, such that one is used for training and the remaining ten are used for testing.  Each of the test images is obtained from the training image of the same class by inverting exactly two pixels.
Since every image pixel is restricted to only one of two colours, each pixel can be intrinsically mapped to a single input node without losing information in the activation transformation.  We note that higher precision inputs, such as a multi-bit grey-scale or colour image, can be easily accommodated by the use of multiple input nodes for each pixel, e.g., eight input nodes may be employed to represent each pixel in an 8-bit image.

In our experiments, network training was performed by minimising the Hamiltonian in Eq.~\eqref{eq:training} and extended ones presented in subsequent sections using an Ising machine.
The scaling factor \(\alpha\) in Eq.~\eqref{eq:training} was fixed to \(\alpha=1\) throughout.
As the Ising machine, we employed a massively-parallel simulated annealing method with Markov chain Monte Carlo optimised for the many-core architecture such as GPU~\cite{okuyama2019}.
The code was executed in single precision on a workstation equipped with an Intel Core i9-9900K CPU (3.60 GHz) and an NVIDIA GeForce RTX 2080 GPU. One annealing calculation was configured with 1,000 independent replicas and evolved for \(N_{\rm MC}=1,000\) Monte Carlo steps for the variable updates. 
The temperature parameter \(T\) in the annealing algorithm was updated after every Monte Carlo step according to the geometric cooling schedule \(T = T_{\rm max} (T_{\rm min}/T_{\rm max})^{t / (N_{\rm MC}-1)}\) for each step \(t=0,\dots,N_{\rm MC}-1\).  The temperature parameters \(T_{\rm min}\) and \(T_{\rm max}\) were tuned by the Nelder-Mead method~\cite{nelder1965}.
This entire annealing procedure was independently repeated for training each BNN.

In Section \ref{sec:architectures}, \ref{sec:sum-of-margins}, and \ref{sec:dropout}, we report on the results of our computational experiments using this setup unless otherwise stated.

\section{Analysis of network architectures}
\label{sec:architectures}

\begin{figure}[htbp]
\centering
\includegraphics[width=.9\linewidth]{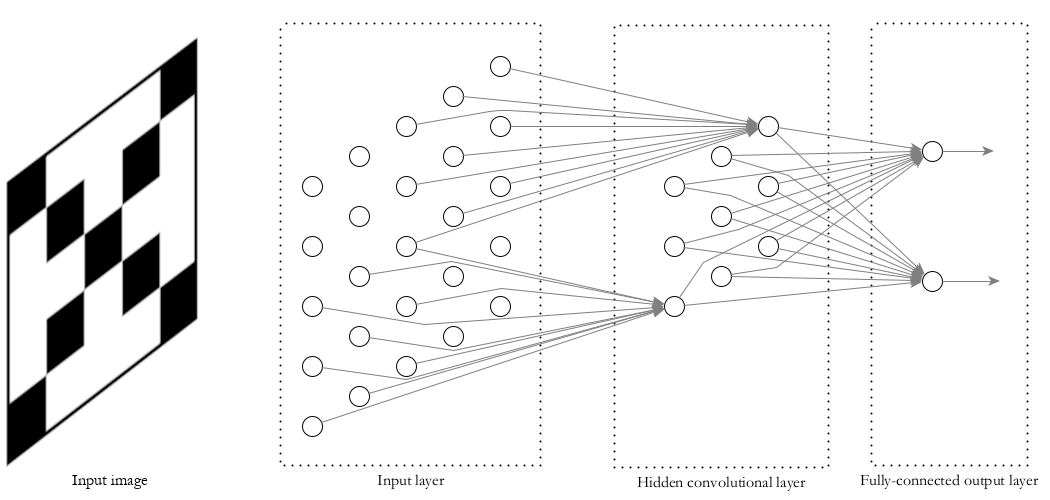}
\caption{\label{fig:conv-network-example}Schematic diagram of Network 2 in Table \ref{tab:networks}, that is a convolutional neural network with \(5\times5\) input nodes, \(3\times3\) convolutional filter giving rise to a \(3\times3\) node hidden convolutional layer, and a 2-node fully-connected output layer.  Note, that for clarity, incoming connections are only shown for two of the nodes in the hidden layer.}
\end{figure}

In this section, the training formulation introduced in previous sections is directly applied to several network architectures. We note that the results provide a baseline for regularisation methods discussed in the latter sections.  Eighteen different network architectures are employed:
eight convolutional networks (having one or two hidden convolutional layers) and ten fully-connected networks (having a single fully-connected hidden layer).  The networks are summarised in Table~\ref{tab:networks}.  A diagram showing the architecture of a convolutional network (Network 2) is shown in Figure~\ref{fig:conv-network-example}.
The convolutional networks apply filters of size \(2\times2\), \(3\times3\), or \(4\times4\) with shared weights.
For all networks, the output layer comprises two neurons, such that any combination of outputs maps to exactly one of the four labels.  E.g., an output of \((-1,-1)\) corresponds to the label \emph{O}, while an output of \((1,-1)\) corresponds to the label \emph{X}.  For each network architecture 200 tests were conducted, each comprising a training run and inference on each of the 40 images in the test data set.  In the following, we present and analyse the details of the experiments described above.

\begin{table}[htbp]
\caption{\label{tab:networks}Summary of network topologies investigated.  Each network architecture is composed of one or two hidden layers, each of which is either convolutional (conv) or fully connected (fc).  fc(\(a\)) indicates a fully-connected layer with \(a\) neurons, while conv(\(a \times b\)) indicates a convolutional layer with an \(a \times b\) filter. conv(\(a \times b \times 2\)) indicates that 2 \(a \times b\) filters are applied in the same layer.  When two hidden layers are used the + operator indicates separation between the two layers.  \emph{Neurons} and \emph{Connections} indicate the number of neurons and connections, respectively, in the network.  \emph{Constraints} indicates the total number of activation constraints \eqref{eq:training-activation-constraint2} and auxiliary substitution constraints \eqref{eq:training-linear-substitution},  \emph{Variables} indicates the total number of binary variables \(v_{ij}\), \(y^k_i\), and \(\psi^k_{ij}\) and integer variables \(\chi^k_j\), respectively, used in the QUBO formulation of the training problem.}
\centering
\resizebox{\linewidth}{!}{
  \begin{tabular}{clccccc}
\toprule
Network & Architecture & Neurons & Connections & \multicolumn{2}{c}{Variables} & Constraints\\[0pt]
\cline{5-6}
        &              & $\vert\N \vert$ & $\vert\A \vert$     & binary & integer              &            \\[0pt]
\midrule
0 & conv(\(2{\times}2\)) & 43 & 96 & 246 & 72 & 200\\[0pt]
1 & conv(\(2{\times}2\)) + fc(4) & 47 & 136 & 466 & 88 & 376\\[0pt]
2 & conv(\(3{\times}3\)) & 36 & 99 & 146 & 44 & 116\\[0pt]
3 & conv(\(3{\times}3{\times}2\)) & 45 & 198 & 290 & 80 & 224\\[0pt]
4 & conv(\(3{\times}3\)) + fc(4) & 40 & 125 & 296 & 60 & 236\\[0pt]
5 & conv(\(4{\times}4\)) & 31 & 72 & 78 & 24 & 56\\[0pt]
6 & conv(\(4{\times}4{\times}2\)) & 35 & 144 & 154 & 40 & 104\\[0pt]
7 & conv(\(4{\times}4\)) + fc(4) & 39 & 168 & 294 & 56 & 216\\[0pt]
\midrule
8 & fc(1) & 28 & 27 & 42 & 12 & 20\\[0pt]
9 & fc(2) & 29 & 54 & 82 & 16 & 32\\[0pt]
10 & fc(3) & 30 & 81 & 122 & 20 & 44\\[0pt]
11 & fc(4) & 31 & 108 & 162 & 24 & 56\\[0pt]
12 & fc(5) & 32 & 135 & 202 & 28 & 68\\[0pt]
13 & fc(6) & 33 & 162 & 242 & 32 & 80\\[0pt]
14 & fc(7) & 34 & 189 & 282 & 36 & 92\\[0pt]
15 & fc(8) & 35 & 216 & 322 & 40 & 104\\[0pt]
16 & fc(9) & 36 & 243 & 362 & 44 & 116\\[0pt]
17 & fc(10) & 37 & 270 & 402 & 48 & 128\\[0pt]
\bottomrule
    \end{tabular}
    }
\end{table}

\begin{figure}[htbp]
\centering
\includegraphics[width=.9\linewidth]{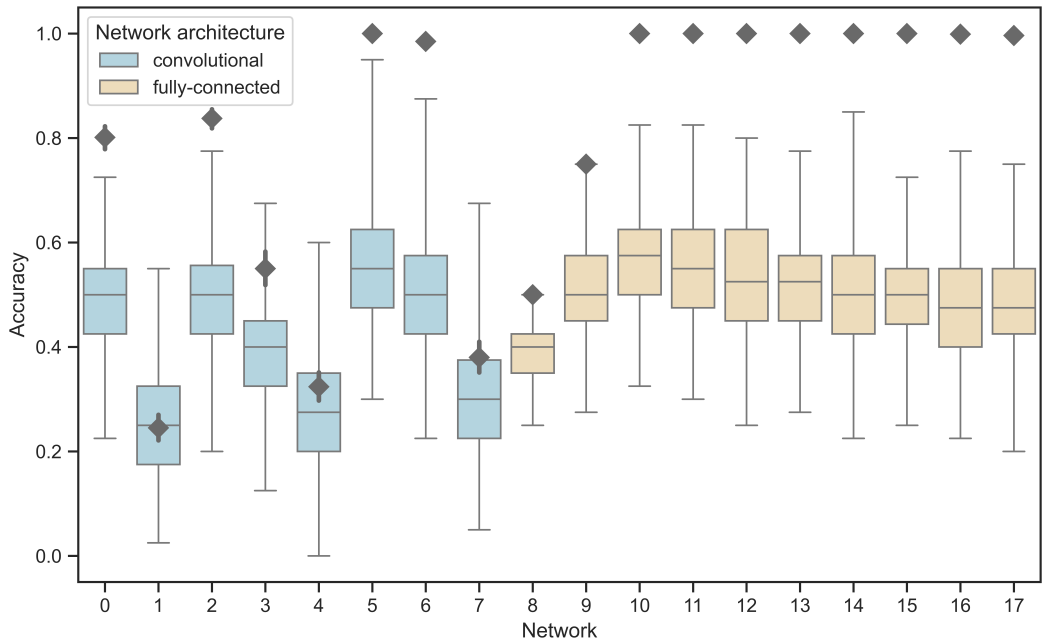}
\caption{\label{fig:boxplot-accuracy-all}Test and training accuracy statistics for the 18 network architectures.  Each box represents the first quartile, median, and third quartile, while whiskers indicate the range of the distribution (minimum to maximum) of the test accuracy.  The grey squares indicate the mean training accuracy, while the attached grey vertical lines indicate the range of training accuracies (minimum to maximum).}
\end{figure}

Table~\ref{tab:results-oneshot} and Figure~\ref{fig:boxplot-accuracy-all} summarise the results by network architecture.  The best performing networks are network 10, 11, and 5, with a mean test accuracy of 57\%, 55\%, and 55\%, respectively.  All three networks have a mean training accuracy of 1.0.
Networks 10 and 11 are both fully-connected networks with three and four neurons, respectively, in the middle hidden layer, while Network 5 is a convolutional network with four neurons in the middle layer.  The best performing networks are all  relatively small networks with a modest number of neurons and connections.  Network 5 is the smallest among the convolutional networks (0-7), while Network 10 and 11 are the smallest networks among the fully-connected networks (8-17), that have at least three neurons in the middle layer.
All the fully-connected networks with at least three nodes in the middle layer obtained a mean training accuracy of or close to 1. 

It is worth remarking, that all networks with optimal solutions to the Hamiltonian, specifically Network 5 and Networks 10–15, satisfy all activation and auxiliary constraints and achieve a training accuracy of 1.0, as expected.
In contrast, the worst performing networks (e.g., Network 1, 3, 4, 7, 8, and 9) all have a large fraction of unsatisfied constraints in the solution of the training problem.  This is a direct cause of low training accuracy, since many of the neuron activations are not correctly applied.  The convolutional networks with poor performance (Network 1, 3, 4, 7) are all relatively large models with more than 200 constraints each.  It is likely that for large models, the low performance stems from the inability of the annealing algorithm to converge on solutions that satisfy all constraints.  On the contrary, the fully-connected networks (Network 8 and 9) are very small networks with few constraints.  In these cases the poor performance is caused by the network's inability to encode the input signals properly.  For instance, Network 8 has only one neuron (and a 1-bit activation) in the middle layer, while two bits are necessary to encode the classification of the input into one of four classes.  Hence, half of the required information is lost in the middle layer resulting in a mean training accuracy of only 0.5.

\begin{figure}
\centering

  \begin{subfigure}{0.45\textwidth}
    \centering
    \includegraphics[width=\linewidth]{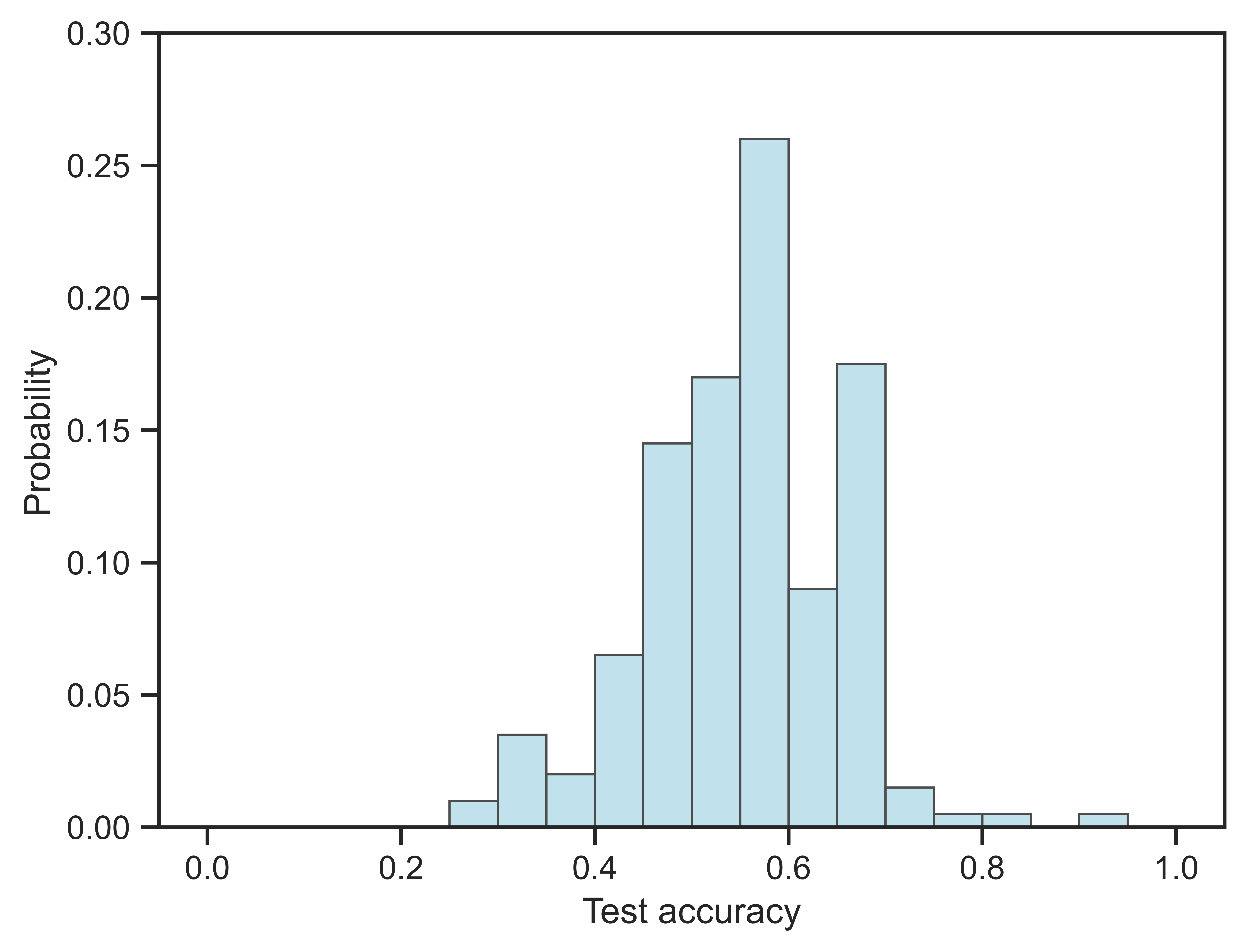}
    \caption{Network 5.}
    \label{fig:hist-net5}
  \end{subfigure}
  \hspace{1cm}
  \begin{subfigure}{0.45\textwidth}
    \centering
    \includegraphics[width=\linewidth]{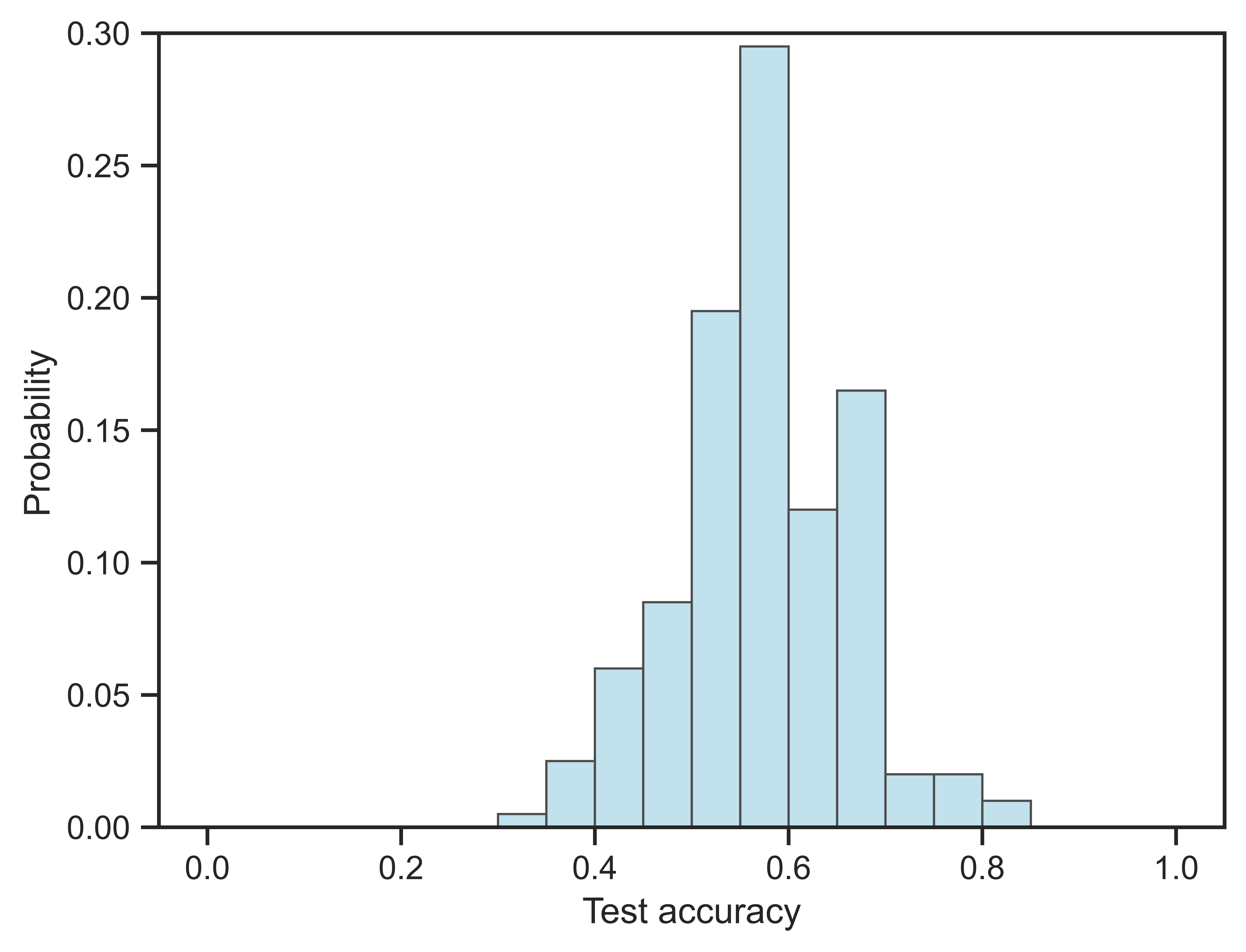}
    \caption{Network 10.}
    \label{fig:hist-net10}
  \end{subfigure}
\caption{Histogram of test accuracy across 200 independent training and inference runs for each of Network~5~\subrefp{fig:hist-net5} and Network~10~\subrefp{fig:hist-net10}.}
\label{fig:hist-architectures}
\end{figure}

Figure~\ref{fig:hist-architectures} shows the histograms of test accuracies obtained from Networks 5 and 10, which represent the best-performing convolutional and fully connected architectures, respectively.
The distributions are uni-modal and concentrated around their mean values of approximately 0.55 for Network 5 and 0.57 for Network 10.  Note that the annealing calculations found optimal solutions for all training runs of both networks but the variation in the test accuracies suggests that, even among optimal solutions, the network parameter configurations of some solutions generalise better than those of others to unseen data.  This observation motivates the discussion in the following sections, where we explore methods for regularisation schemes for QUBO-based BNN to encourage solutions that achieve better generalisation performance.

\begin{table}[htbp]
\caption{\label{tab:results-oneshot}Summary statistics for test and training accuracy for each of the 18 network architectures investigated. Each entry represents results obtained from 200 independent training and inference runs on the four-class image classification task.  The table lists the minimum, maximum, mean, and median test accuracies, as well as the corresponding mean training accuracies and the average fraction of unsatisfied constraints \eqref{eq:training-activation-constraint2} and \eqref{eq:training-linear-substitution}.  The entries highlighted in bold, Networks 5 and 10, indicate the best test accuracy among the convolutional and fully connected architectures, respectively.}
\centering
\resizebox{\linewidth}{!}{
\begin{tabular}{ccccccc}
\toprule
Network & \multicolumn{4}{c}{Test accuracy}           & Training accuracy & Unsatisfied constraints\\[0pt]
\cline{2-5}
        &  min.     & max.      & mean      & median  & mean              & mean fraction ($\%$)\\[0pt]
\midrule
0 & 0.225 & 0.725 & 0.487 & 0.500 & 0.801 & 0.6875\\[0pt]
1 & 0.025 & 0.550 & 0.254 & 0.250 & 0.245 & 4.6782\\[0pt]
2 & 0.200 & 0.775 & 0.495 & 0.500 & 0.838 & 0.7845\\[0pt]
3 & 0.125 & 0.675 & 0.390 & 0.400 & 0.550 & 2.8549\\[0pt]
4 & 0.000 & 0.600 & 0.274 & 0.275 & 0.324 & 2.6483\\[0pt]
\bf5 & \bf0.300 & \bf0.950 & \bf0.550 & \bf0.550 & \bf1.000 & \bf 0.0000\\[0pt]
6 & 0.225 & 0.875 & 0.505 & 0.500 & 0.985 & 0.0865\\[0pt]
7 & 0.050 & 0.675 & 0.301 & 0.300 & 0.380 & 2.2917\\[0pt]
\midrule
8 & 0.250 & 0.500 & 0.392 & 0.400 & 0.500 & 10.0000\\[0pt]
9 & 0.275 & 0.750 & 0.509 & 0.500 & 0.750 & 3.1250\\[0pt]
\bf 10 & \bf 0.325 & \bf 0.825 & \bf 0.569 & \bf 0.575 & \bf 1.000 & \bf 0.0000\\[0pt]
11 & 0.300 & 0.825 & 0.552 & 0.550 & 1.000 & 0.0000\\[0pt]
12 & 0.250 & 0.800 & 0.533 & 0.525 & 1.000 & 0.0000\\[0pt]
13 & 0.275 & 0.775 & 0.512 & 0.525 & 1.000 & 0.0000\\[0pt]
14 & 0.225 & 0.850 & 0.505 & 0.500 & 1.000 & 0.0000\\[0pt]
15 & 0.250 & 0.725 & 0.500 & 0.500 & 1.000 & 0.0000\\[0pt]
16 & 0.225 & 0.775 & 0.482 & 0.475 & 0.999 & 0.0086\\[0pt]
17 & 0.200 & 0.750 & 0.486 & 0.475 & 0.996 & 0.0273\\[0pt]
\bottomrule
\end{tabular}
}
\end{table}

\section{Maximising sum of margins}
\label{sec:sum-of-margins}

In classification problems, hyper planes define decision boundaries between classes and determine how an input is classified.  Margins define the minimum distance from any data point to the decision boundary, and thus reflect the model's robustness to input perturbations. Extending this idea to neural networks, the neuron margin quantifies how much the pre-activation input must change to flip the neuron's output.
As maximising neuron margins has been shown to enhance generalisation and robustness, efforts have been undertaken to incorporate margin-based loss functions in the training of neural networks~\cite{sun2015,elsayed2018,jiang2018,toro2019}.  Building upon these observations, we introduce a QUBO-based formulation that promotes neuron margins in BNNs.

\subsection{Preliminary analysis and model formulation}
Following~\cite{toro2019}, the neuron margin for neuron \(j \in \N\) is defined as the minimum absolute value of its pre-activation across all training samples \(k \in \K\).

\begin{align}
  \min_{k\in \K} &\left|  b_j+ \sum_{i\in P_j} w_{ij} x^k_i \right|  \label{eq:margin-def}
\end{align}

This value measures the smallest perturbation required to change the neuron's activation.
Summing the margins over all non-input nodes, we obtain a network level indicator of the \emph{robustness} of the network given the training data,

\begin{align}
S_1 &= \sum_{j\in \N\setminus \N_0} \min_{k\in \K} \left|  b_j+ \sum_{i\in P_j} w_{ij} x^k_i \right|.  \label{eq:sumofmargins}
\end{align}

However, \(S_1\) is not amenable to optimisation using QUBO expressions as it involves nested minimisation.  To obtain a tractable expression, we define the alternative proxy indicator \(S_2\), where instead of  taking the minimum over all data points, we take the summation as follows,

\begin{align}
S_2 &= \sum_{j\in \N\setminus \N_0} \sum_{k\in \K} \left|  b_j+ \sum_{i\in P_j} w_{ij} x^k_i \right|.  \label{eq:sumofmargins2}
\end{align}

\begin{figure}
  \centering

  \begin{subfigure}{0.45\textwidth}
    \centering
    \includegraphics[width=\linewidth]{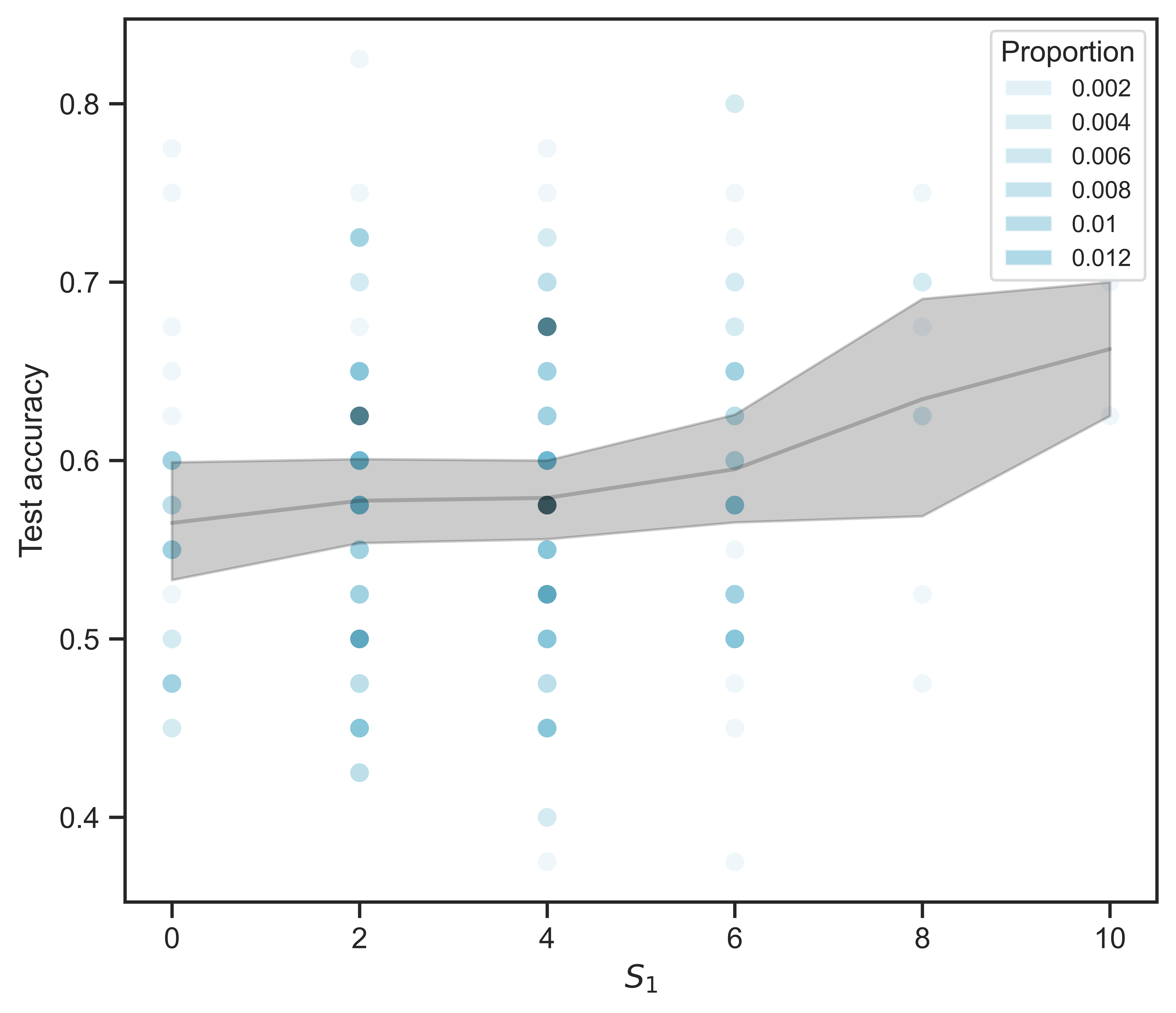}
    \caption{ }
    \label{fig:som_S1_gamma0}
  \end{subfigure}
  \hspace{1cm}
  \begin{subfigure}{0.45\textwidth}
    \centering
    \includegraphics[width=\linewidth]{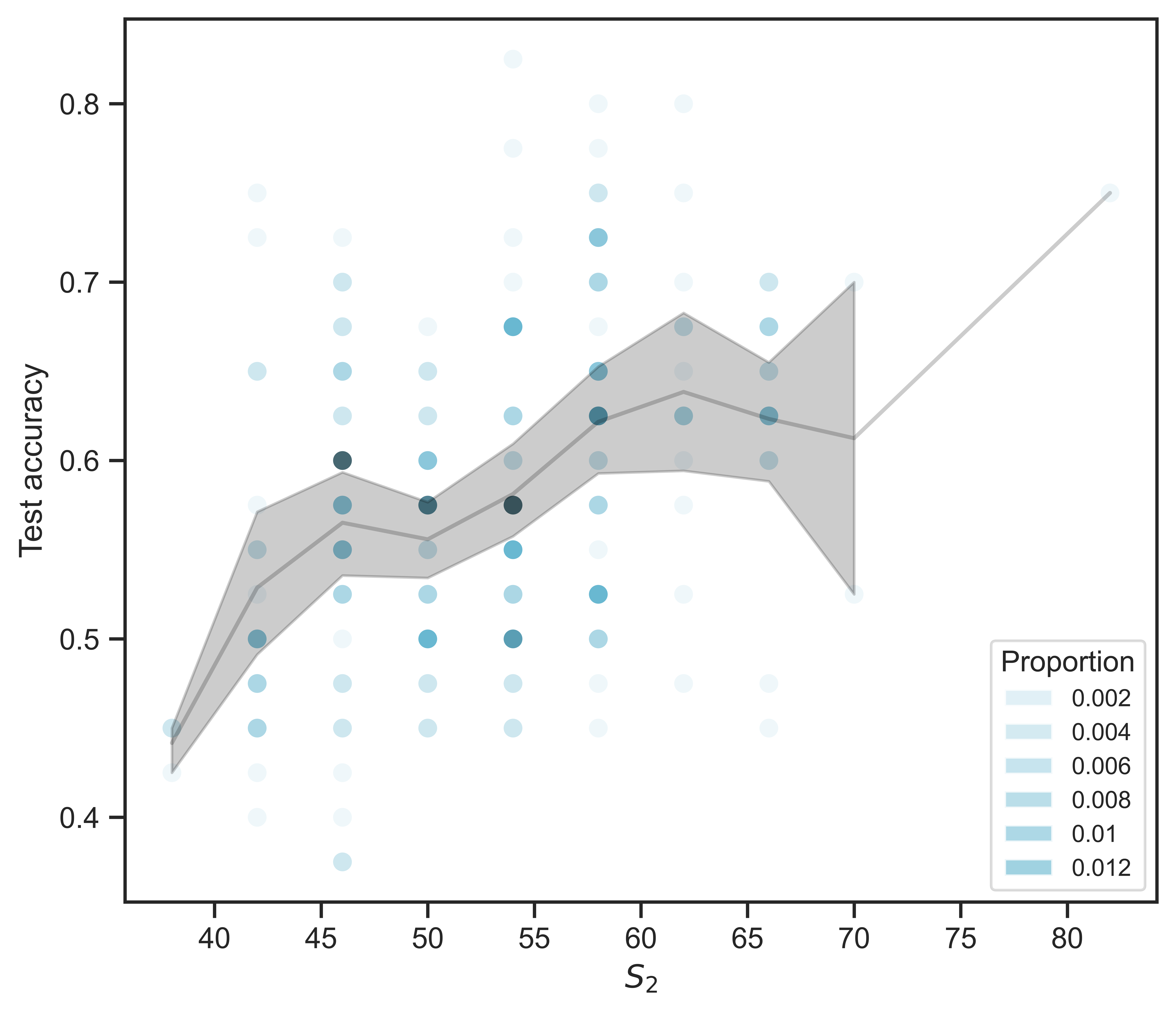}
    \caption{ }
    \label{fig:som_S2_gamma0}
  \end{subfigure}

  \caption{Test accuracy vs.~sum of margins $S_1$~\subrefp{fig:som_S1_gamma0} and $S_2$~\subrefp{fig:som_S2_gamma0} evaluated across 200 independent training and inference runs using Network 10 without regularisation.  Dots indicate observations, while the grey line indicates the mean accuracy and the shaded band represents the 95\% confidence interval computed using the standard error of the mean.  The lightness of the dots indicates the fraction of all 200 tests that have the same combination of test accuracy and $S_1$-, respectively, $S_2$-value.}
  \label{fig:som}
\end{figure}

In Figure~\ref{fig:som}, we plot the test accuracy against \(S_1\) and \(S_2\) for Network 10 trained without regularisation (as in Section~\ref{sec:architectures}) for 200 runs of the algorithm. This gives us an indication of how well the two indicators perform as predictors for generalisation.  Since, \(S_1\) relies on taking the minimum of the absolute value of the pre-activation values, the number of possible outcomes is relatively limited and in our case is restricted to one of six values (0, 2, 4, 6, 8, and 10), while the domain of \(S_2\) on the other hand is more granular.  Furthermore, from Figure~\ref{fig:som}, \(S_2\) appears to exhibit a clearly positive correlation with test accuracy---that is, a higher value of \(S_2\)  seems to indicate a higher accuracy on unseen data for this experiment.  In the following, we show how to formulate the proxy margin indicator \(S_2\) as a QUBO expression and incorporate it into the training Hamiltonian.

We note first that for every input data point \(k\in\K\) and non-input node \(j\in \N\setminus \N_0\), we have from the activation constraint \eqref{eq:training-activation-constraint2} and the definition of \(\rho_j\),

\begin{align*}
\frac{b_j+1}{2}+ \sum_{i\in P_j} \frac{w_{ij} x^k_i+1}{2} = 2^{n_j} y_j^k + \chi^k_j - \const
\end{align*}

Therefore, the pre-activation value for a non-input node \(j\) is given by 

\begin{align}
b_j+ \sum_{i\in P_j} w_{ij} x^k_i &= 2 \left(2^{n_j} y_j^k + \chi^k_j - \const \right) - |P_j| - 1  \label{eq:som1}
\end{align}

Since, \(y_j=1\) when the pre-activation value is strictly positive and \(y_j=0\) otherwise, we have that the absolute value of the pre-activation of node \(j\in \N \setminus \N_0\) for data point \(k\) can be written as

\begin{align}
  \left| b_j+ \sum_{i\in P_j} w_{ij} x^k_i \right|    &=  (2y_j^k-1) \left( b_j+ \sum_{i\in P_j} w_{ij} x^k_i \right) \label{eq:som2}
\end{align}

From Eqs.~\eqref{eq:som1} and \eqref{eq:som2}, we can re-write the absolute value of the pre-activation as,

\begin{align}
  \left| b_j+ \sum_{i\in P_j} w_{ij} x^k_i \right| &= (2y_j^k-1) \left(2 \left(2^n y_j^k + \chi^k_j - \const \right)-|P_j| - 1 \right)  \label{eq:som-expression}
\end{align}

The right hand side of Eq.~\eqref{eq:som-expression} is quadratic in \(y_j\) and can readily be incorporated in the QUBO formulation.  We therefore introduce the following QUBO penalty term,

\begin{align}
      H_{som}    &=  \sum_{j\in \N\setminus \N_0} \sum_{k\in \K} (2y^k_j-1) \left(2 \left(2^n y^k_j + \chi^k_j - \const \right)-|P_j| - 1 \right)  \label{eq:som-hamiltonian}
\end{align}

To maximise \(S_2\), we subtract \(H_{som}\) from the original Hamiltonian \(H_{train}\) to obtain the augmented energy function,

\begin{align}
H^*_{train}= H_1 + \alpha H_2 - \gamma H_{som} \label{eq:training-som}
\end{align}

where \(\alpha\) and \(\gamma\) are scaling factors.

\subsection{Computational results}
\begin{figure}
  \centering

  \begin{subfigure}{0.45\textwidth}
    \centering
    \includegraphics[width=\linewidth]{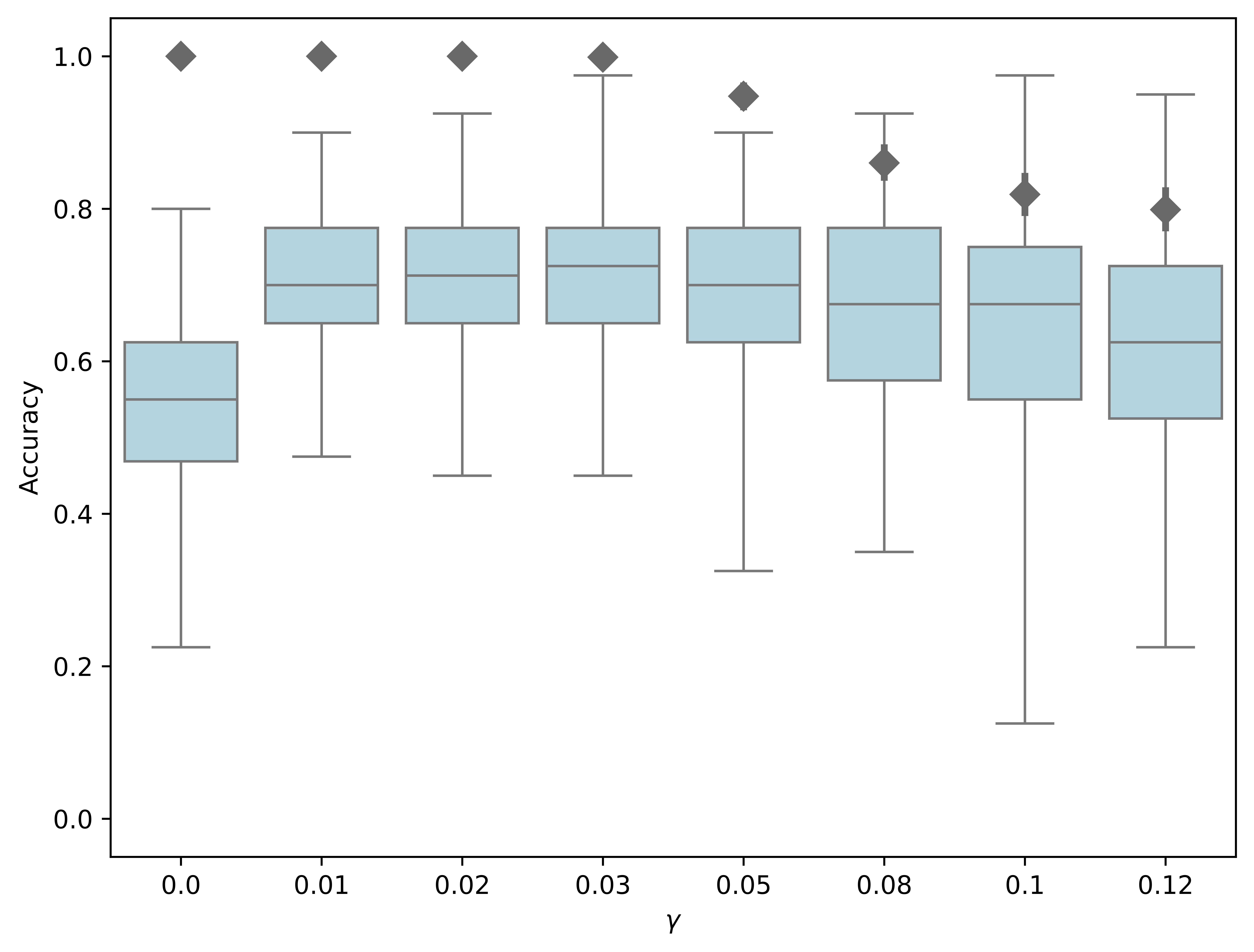}
    \caption{Network 5.}
    \label{fig:som-results-network5}
  \end{subfigure}
  \hspace{1cm}
  \begin{subfigure}{0.45\textwidth}
    \centering
    \includegraphics[width=\linewidth]{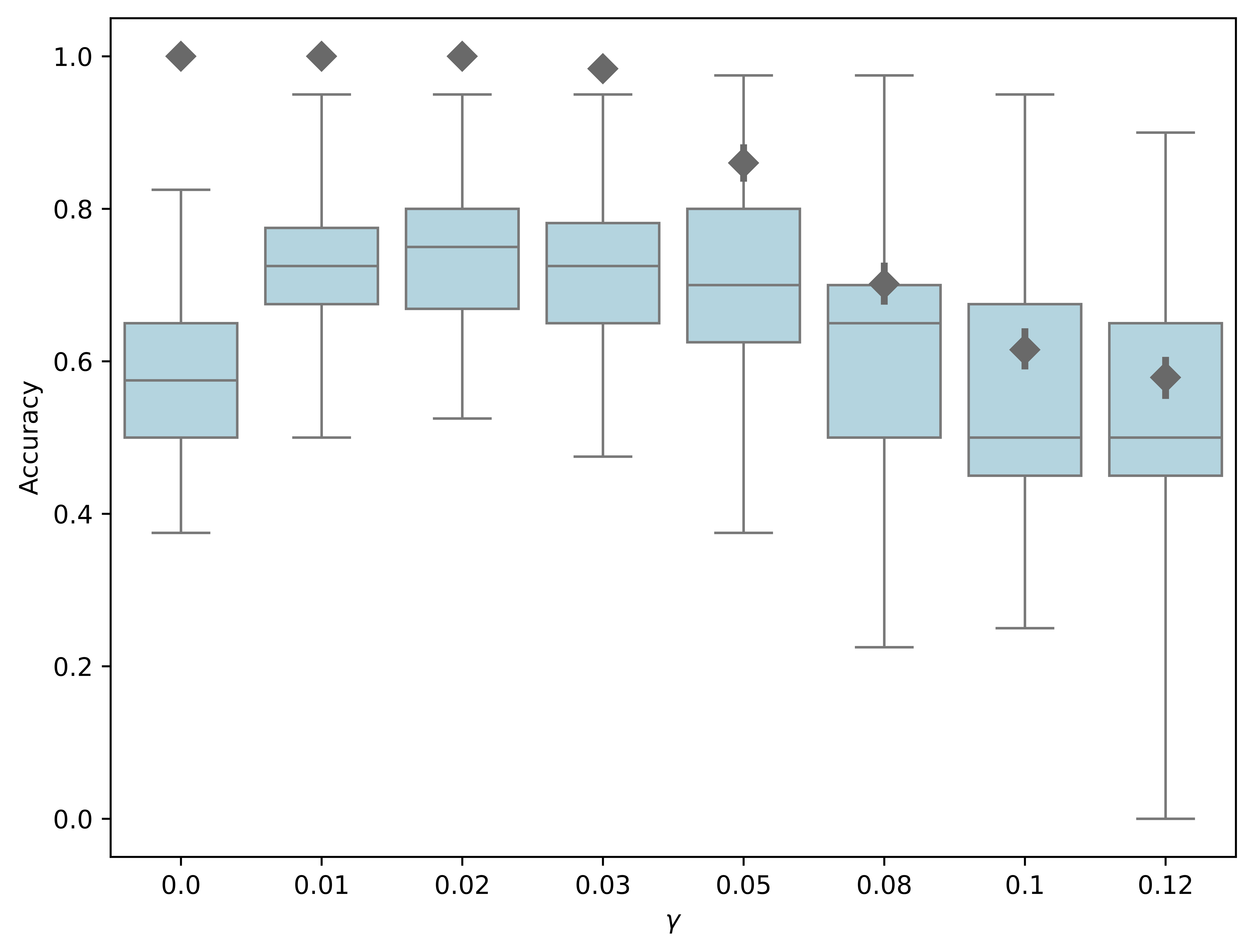}
    \caption{Network 10.}
    \label{fig:som-results-network10}
  \end{subfigure}

  \caption{Distribution of test accuracy for Network~5~\subrefp{fig:som-results-network5} and Network~10~\subrefp{fig:som-results-network10}, evaluated across 200 training and inference runs, for each value of the margin-scaling factor $\gamma$.  Each box represents the first quartile, median, and third quartile, while the whiskers indicate the  minimum and the maximum of the distribution. Grey squares and the corresponding vertical lines denote the mean and the range (minimum to maximum) of the training accuracy for the corresponding $\gamma$, respectively.}
  \label{fig:som-results}
\end{figure}

To examine the effect of maximising \(S_2\) on QUBO-based BNNs, we conduct computational experiments on Network 5 and Network 10.  The experiments are conducted as in Section~\ref{sec:architectures} with the objective function replaced by \(H^*_{train}\) \eqref{eq:training-som}.
In all experiments, \(\alpha=1\) is kept constant, while varying the value of the scaling factor \(\gamma\).  For each value of  \(\gamma\), we perform 200 runs of training and inference.  The results are summarised in Table~\ref{tab:som-summary} and Figure~\ref{fig:som-results}.  Experiments with \(\gamma=0\) are equivalent to the experiments performed in Section~\ref{sec:architectures} without regularisation and serve as a baseline for comparison.
Table~\ref{tab:som-summary} and Figure~\ref{fig:som-results} indicate that the performance of the BNNs on classification tasks improves significantly when introducing margin-based regularisation maximising \(S_2\).
For Network 5, the best results are achieved with \(\gamma=0.03\), which increases the mean accuracy from 54\% to 71\% corresponding to a performance improvement of 31\%.  Similarly, for Network 10, the highest mean accuracy is achieved with \(\gamma=0.02\) leading to an increase in mean accuracy from 58\% to 73\% or a 26\% relative increase in accuracy compared to the baseline without regularisation.

\begin{figure}
  \centering

    \begin{subfigure}{0.45\textwidth}
      \centering
      \includegraphics[width=\linewidth]{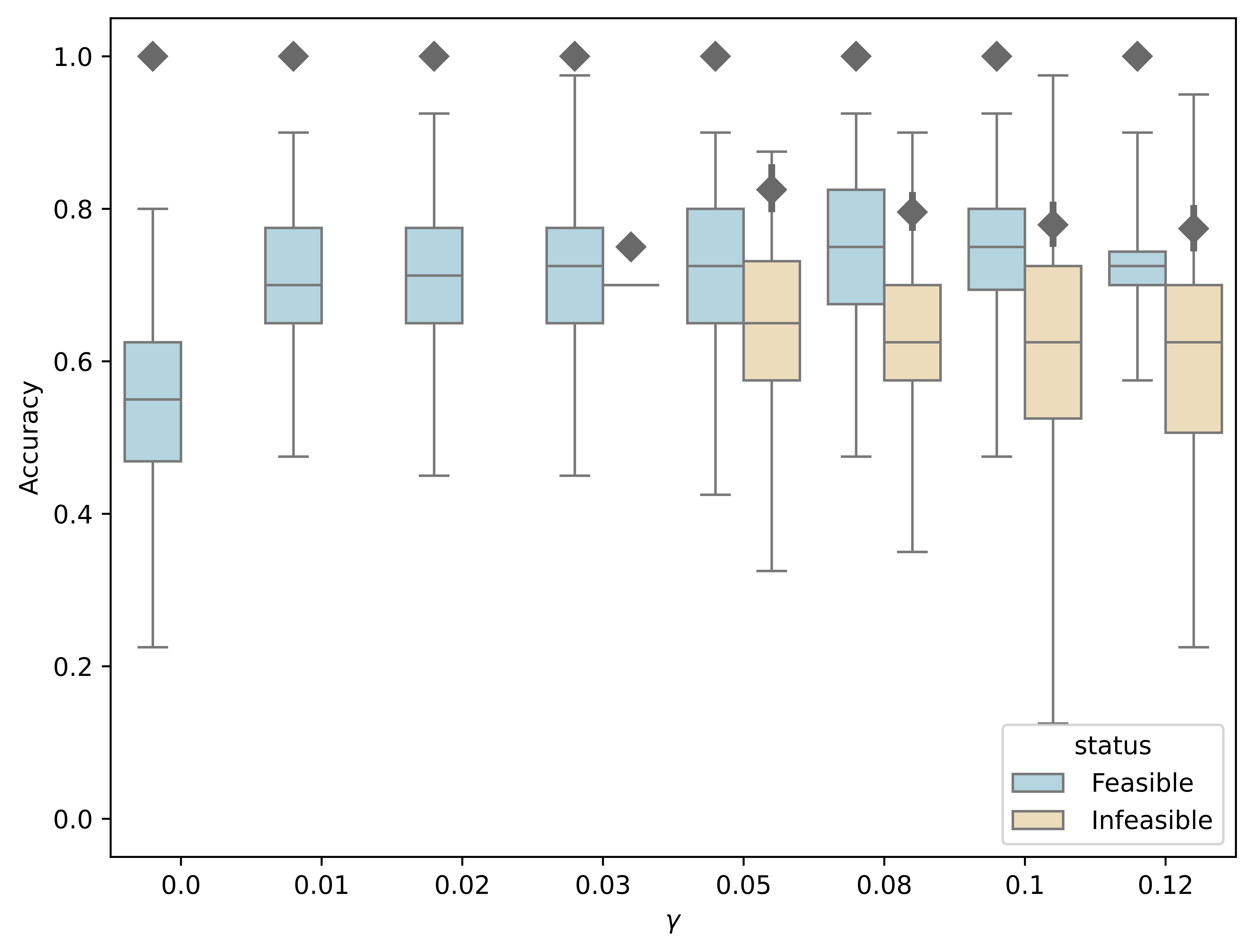}
      \caption{Network 5.}
    \label{fig:som-results-grouped-network5}
  \end{subfigure}
  \hspace{1cm}
  \begin{subfigure}{0.45\textwidth}
    \centering
    \includegraphics[width=\linewidth]{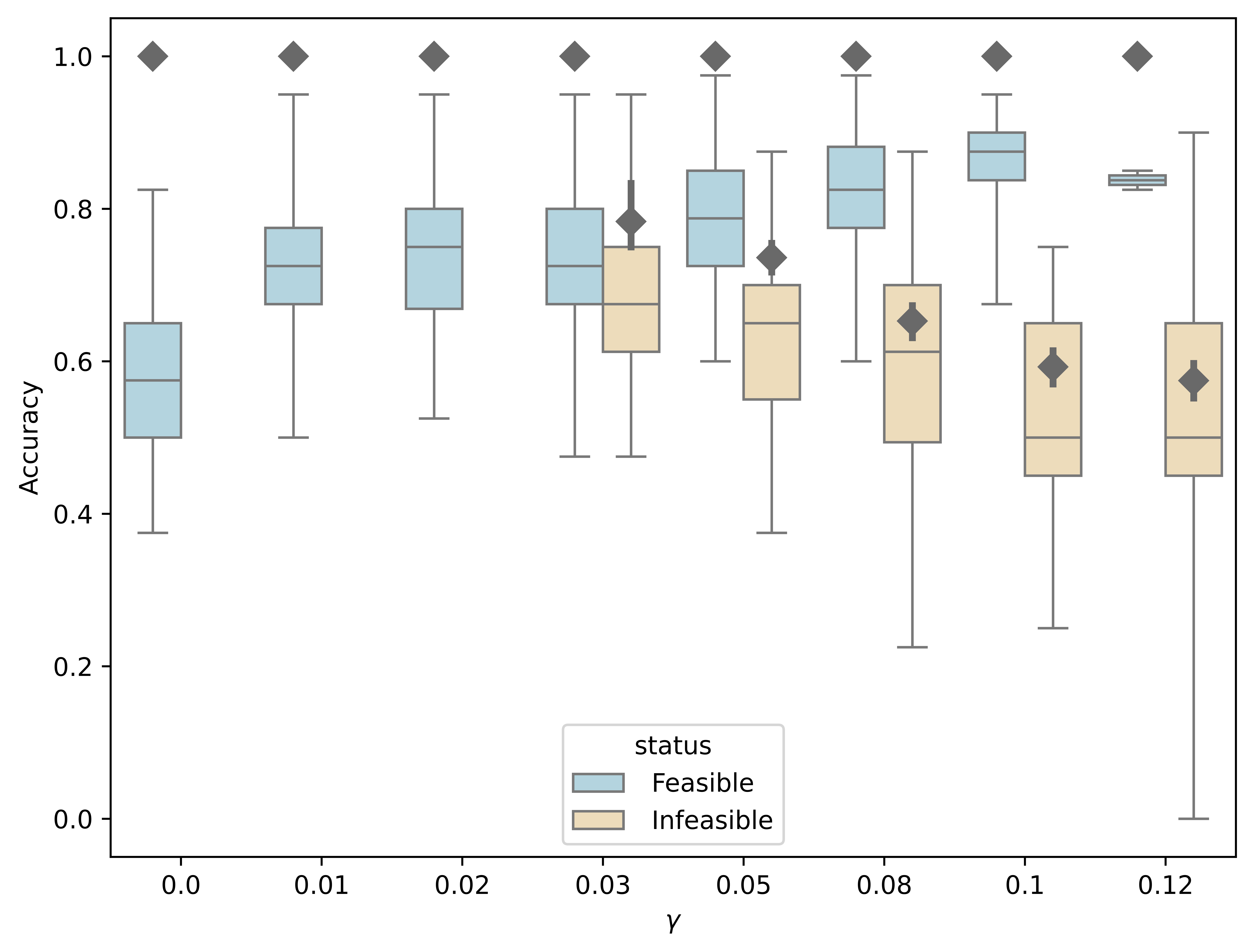}
    \caption{Network 10.}
    \label{fig:som-results-grouped-network10}
  \end{subfigure}

\caption{Distribution of test accuracy for Network~5~\subrefp{fig:som-results-grouped-network5} and Network~10~\subrefp{fig:som-results-grouped-network5}, evaluated across 200 independent training and inference runs, grouped by feasibility (status) of the solutions obtained, for each value of the margin-scaling factor $\gamma$.  Each box represents the first quartile, median, and third quartile, while the whiskers indicate the  minimum and the maximum of the distribution. Grey squares and corresponding vertical lines denote the mean and the range (minimum to maximum) of training accuracy for the corresponding $\gamma$, respectively.}
\label{fig:som-results-grouped}
\end{figure}

When the scaling factor \(\gamma\) becomes too large, the performance tends to decrease.  As shown in Table~\ref{tab:som-summary}, for \(\gamma\geq0.03\) the annealing algorithm sometimes produces infeasible solutions that do not satisfy either activation or auxiliary constraints (i.e., \(H_1 > 0\) or \(H_2 > 0\)).  Furthermore, the proportion of infeasible solutions increases with increasing \(\gamma\). 
This behavior stems from the dominant contribution of \(H_{\mathrm{som}}\) in the Hamiltonian, which emphasises margin maximisation rather than constraint satisfaction.  As infeasible solutions may violate the activation constraints, the networks derived from such solutions exhibit low training accuracy, thereby lowering the overall average performance across all trials.
However, when only feasible solutions are considered, the test accuracy increases monotonically up to approximately \(\gamma \approx 0.1\) (Figure~\ref{fig:som-results-grouped}).  In particular, for Network~10, the mean accuracy reaches almost 90\%.

\begin{figure}
\centering
\includegraphics[width=0.75\linewidth]{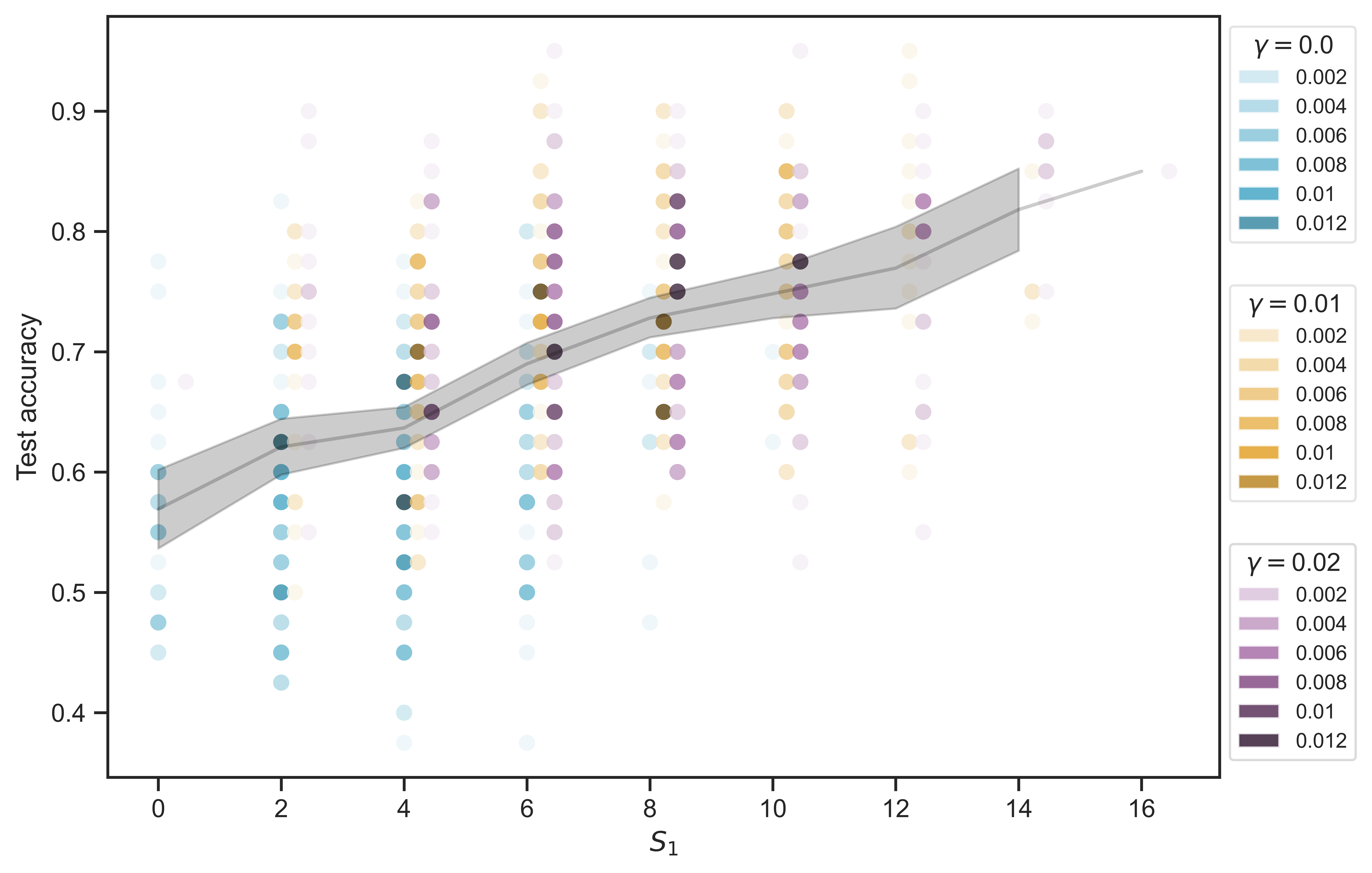}
\caption{Test accuracy vs.~sum of margins $S_1$ for all tests using Network 10 with $\gamma=0.00, 0.01, 0.02$.  Dots indicate observations, while the grey line (band) shows the mean value (95\% confidence interval) of the observations.  The lightness of the dots indicate for each $\gamma$-value, the fraction of tests that have the same combination of test accuracy and $S_1$-value. For clarity, plots for $\gamma\geq0.01$ are shifted slightly to the right to avoid overlap.}
\label{fig:acc-vs-som1}
\end{figure}

\begin{figure}
\centering
\includegraphics[width=0.75\linewidth]{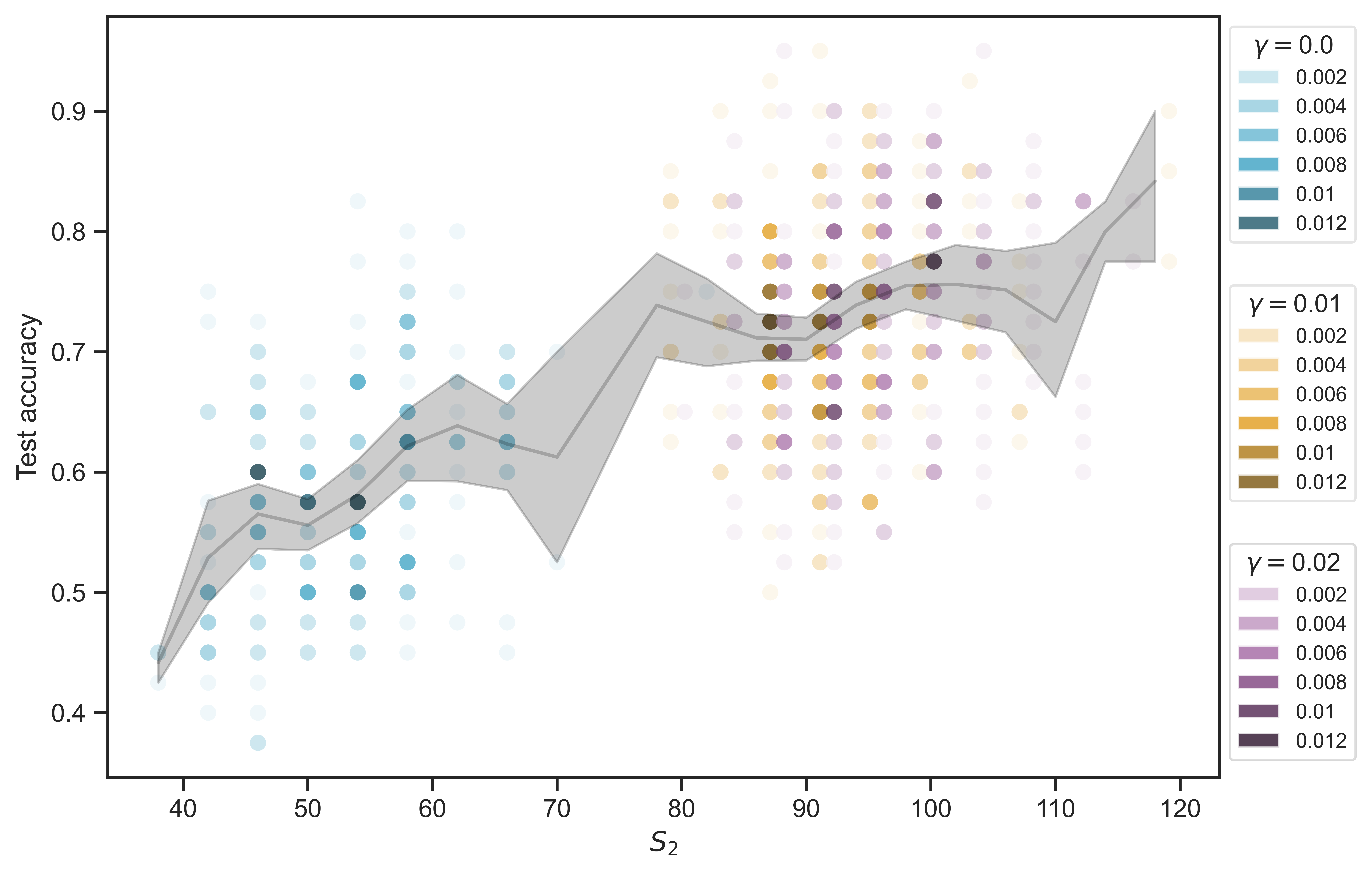}
\caption{Test accuracy vs.~sum of margins $S_2$ for all tests using Network 10 with $\gamma=0.00, 0.01, 0.02$.  Dots indicate observations, while the grey line (band) shows the mean value (95\% confidence interval) of the observations.  The lightness of the dots indicate for each $\gamma$-value, the fraction of tests that have the same combination of test accuracy and $S_2$-value. For clarity, plots for $\gamma\geq0.01$ are shifted slightly to the right to avoid overlap.}
\label{fig:acc-vs-som2}
\end{figure}

These results indicate that maximising margins yields solutions with improved generalisation capability.
Figure~\ref{fig:acc-vs-som1} and Figure~\ref{fig:acc-vs-som2} show the test accuracy and the margin indicators $S_1$ and $S_2$, respectively, for Network~10 with \(\gamma\) = 0.00, 0.01,  and 0.02, for which only feasible solutions are obtained by the annealing runs.  Although the results with \(\gamma=0\) in Figure~\ref{fig:som_S1_gamma0} did not show a clear correlation between the test accuracy and \(S_1\),  Figure~\ref{fig:acc-vs-som1} shows that including the results obtained with $\gamma = 0.01$ and $0.02$ results in a clear positive correlation.  This is consistent with previous studies regarding the usefulness of \(S_1\) as a regularisation measure.  Meanwhile, as illustrated in Figure~\ref{fig:acc-vs-som2}, the correspondingly strong positive correlation with \(S_2\) indicates that \(S_2\) serves as an effective approximation in the QUBO formulation.  This might, however, be attributed to the characteristics of the dataset and the small size of the models and further research is necessary to confirm the effectiveness on other data instances.

\begin{table}[htbp]
\caption{\label{tab:som-summary}Summary statistics for test and training accuracy for different sum of margin factors \(\gamma\) for Network~5 (convolutional) and Network~10 (fully connected). Each entry represents results obtained from 200 independent training and inference runs on the four-class image classification task.  The table lists the minimum, maximum, mean, and median test accuracies, the corresponding mean training accuracies, the mean fraction of unsatisfied constraints \eqref{eq:training-linear-substitution} and \eqref{eq:training-activation-constraint2}, as well as, the mean $S_1$- and $S_2$-values.}
\centering
\resizebox{\linewidth}{!}{
\begin{tabular}{cccccccccc}
\toprule
Network & \(\gamma\) & \multicolumn{4}{c}{Test accuracy} & Training acc. & Unsatisfied constr. & \(S_1\) & \(S_2\)\\[0pt]
\cline{3-6}
        &       &  min. & max. & mean & median & mean & mean fraction (\%) & mean & mean\\[0pt]
\midrule
5 & 0.00 & 0.225 & 0.800 & 0.540 & 0.550 & 1.000 & 0.00 & 7.05 & 59.98\\[0pt]
5 & 0.01 & 0.475 & 0.900 & 0.706 & 0.700 & 1.000 & 0.00 & 10.04 & 88.48\\[0pt]
5 & 0.02 & 0.450 & 0.925 & 0.704 & 0.713 & 1.000 & 0.00 & 10.34 & 89.81\\[0pt]
5 & 0.03 & 0.450 & 0.975 & 0.714 & 0.725 & 0.999 & 0.01 & 9.91 & 89.57\\[0pt]
5 & 0.05 & 0.325 & 0.900 & 0.694 & 0.700 & 0.948 & 0.54 & 10.84 & 95.60\\[0pt]
5 & 0.08 & 0.350 & 0.925 & 0.673 & 0.675 & 0.860 & 1.40 & 11.45 & 102.67\\[0pt]
5 & 0.10 & 0.125 & 0.975 & 0.644 & 0.675 & 0.819 & 2.12 & 11.32 & 106.46\\[0pt]
5 & 0.12 & 0.225 & 0.950 & 0.620 & 0.625 & 0.799 & 2.58 & 10.94 & 107.66\\[0pt]
\midrule
10 & 0.00 & 0.375 & 0.825 & 0.583 & 0.575 & 1.000 & 0.00 & 3.47 & 52.92\\[0pt]
10 & 0.01 & 0.500 & 0.950 & 0.728 & 0.725 & 1.000 & 0.00 & 6.86 & 91.32\\[0pt]
10 & 0.02 & 0.525 & 0.950 & 0.734 & 0.750 & 1.000 & 0.00 & 7.45 & 93.86\\[0pt]
10 & 0.03 & 0.475 & 0.950 & 0.724 & 0.725 & 0.984 & 0.17 & 7.54 & 96.52\\[0pt]
10 & 0.05 & 0.375 & 0.975 & 0.703 & 0.700 & 0.860 & 1.39 & 9.93 & 112.06\\[0pt]
10 & 0.08 & 0.225 & 0.975 & 0.617 & 0.650 & 0.701 & 3.11 & 12.50 & 127.94\\[0pt]
10 & 0.10 & 0.250 & 0.950 & 0.557 & 0.500 & 0.615 & 4.38 & 14.29 & 137.02\\[0pt]
10 & 0.12 & 0.000 & 0.900 & 0.532 & 0.500 & 0.579 & 5.57 & 14.79 & 141.66\\[0pt]
\bottomrule
\end{tabular}
}
\end{table}

\section{A dropout algorithm}
\label{sec:dropout}

Motivated by improvements by the regularisation method in Section \ref{sec:sum-of-margins}, we further investigate an alternative regularisation method.
The approach taken in Sections~\ref{sec:architectures} and \ref{sec:sum-of-margins} relies
on training the neural network by solving a single optimisation problem that includes all training data at once.  When the size of the training data set is large, finding good quality solutions in a single annealing run becomes intractable as the sizes of the models~\eqref{eq:training} and \eqref{eq:training-som} increase with the size of the training data set. 
To mitigate this scalability issue, for example, Georgiev~\cite{georgiev2023} proposed an iterative algorithm that gradually updates network parameters indirectly by iteratively learning from batches of input data. 
In this work, we propose a dropout-like regularisation method that randomly removes neurons and utilises batched training data, thereby reducing the problem size while maintaining regularisation benefits.  We note that this approach dynamically alters the network structure, leveraging the flexibility of our QUBO formulation supporting arbitrary network topologies.

\subsection{Model and algorithm}
In the following, we propose an iterative algorithm for gradually updating network parameters by solving several QUBO problems in sequence.  To this end, we augment the training Hamiltonian to allow for additional (linear) penalty terms on the weights and biases.  We define the augmented Hamiltonian by,

\begin{align}
\minimise H^+_{train}  =  H_1 + \alpha H_2 - H_{ext} \label{eq:augmented-training}
\end{align}

where \(H_{ext} = \sum_{j\in\N} \left( c^b_j b_j+ \sum_{i\in P_j} c^w_{ij}w_{ij} \right)\).

The Hamiltonian \(H_{ext}\) acts as an external weight bias term such that large positive values of \(c^b_j\) for some neuron \(j\) will tend to favour a bias value of \(b_j=1\), whereas a large negative value will favour \(b_j=-1\) and similarly for the weight parameters. 
The algorithm relies on a vector of weight factors \(c^w\) and a vector of bias factors \(c^b\), which are both initialised to the $0$-vectors resulting in \(H_{ext}=0\) initially.
Note, that the scaling of \(H_{ext}\) relative to \(H_1\) and \(H_2\) is implicitly defined by the coefficient vectors \(c^b\) and \(c^w\).

In each iteration of the algorithm, a new configuration \(G'=(\N', \A')\) of the network is obtained by removing a randomly selected subset of nodes \(\N_{drop}\) from the original network \(G\), such that \(\N'=\N \setminus \N_{drop}\) and \(\A'= \{(i,j) \in \A | i,j \notin \N_{drop} \}\).
The new network is then trained on some batch of the training data by minimising the function \(H^+_{train}\) defined on the reduced network \(G'\).
Note, that the removal of nodes from the network affects the number of variables as we are no longer calculating the bias and activations for removed nodes, as well as the weights of their incident connections.

The network parameters (weights and biases) obtained from this optimisation guide the update of the weight and bias factors \(c^w\) and \(c^b\) associated with the reduced network, such that for each neuron \(j\in \N'\), \(c^b_j \leftarrow g(c^b_{j}, b_{j})\) and for each connection \((i,j)\in\A'\), \(c^w_{ij} \leftarrow g(c^w_{ij}, w_{ij})\), for some update function \(g\).  In this paper, we assume the following update function,

\begin{align}
g(a,b) = a + \learnrate \beta^{n^{USC}}b
\end{align}

where \(\learnrate\) is a constant learning factor, \(\beta\leq1\) is a constant factor, and \(n^{USC}\) is the number of unsatisfied constraints in the current solution.
This update function is designed to balance adaptability and stability.  The learning rate \(\learnrate\) controls the step size of the update, determining how strongly the new parameter values influence the bias and weight factors.
The exponential term \(\beta^{n^{USC}}\) acts as an adaptive weight, which exponentially suppresses the update magnitude when constraints are unsatisfied.  Thus, the bias and weight factors increase (or decrease) proportionally to the learning rate and decay exponentially with the number of unsatisfied constraints.
Note that this update rule is a generalised version of the rule presented in \cite{georgiev2023} that explicitly considers the distance from the current solution to a feasible one, and that the two update rules are equivalent for \(\beta=1\).
The number of iterations \(M\) is set by the user.   After the final iteration the full network \(\N\) is trained without removing any nodes.  The algorithm is outlined in Algorithm~\ref{alg:dropout}. 

\begin{algorithm}
    \caption{A dropout algorithm for QUBO-based neural networks.}
    \label{alg:dropout}
    \begin{algorithmic}
    \State $c^w \gets 0$  \Comment{initialise weight factors}
    \State $c^b \gets 0$  \Comment{initialise bias factors}
    \For{$i \gets 1$ to $M$}
    \State $\K_i \gets$  \Call{GetTrainingBatch}{$i$}  \Comment{get training data batch $i$}
    \State $G' \gets$ \Call{ModifyNetwork}{$G$} \Comment{modify network, e.g., by node removal}
    \State $w,b \gets \min_{G',\K_i} H^+_{train}(c^w, c^b)$ \Comment{train network} 
    \State {$c^w \gets g(c^w, w)$} \Comment{update weight factors}
    \State {$c^b \gets g(c^b, b)$} \Comment{update bias factors}
    \EndFor
    \State $w,b \gets \min_{\N,\K_{M+1}} H^+_{train}(c^w, c^b)$ \Comment{train final batch without dropout} 
    \end{algorithmic}
\end{algorithm}

\subsection{Computational results}

We perform computational experiments on the data set described in Section~\ref{sec:setup} and depicted in Figure~\ref{fig:dataset}.  As in the previous experiments, the training data is limited to the first four images, while testing is performed on the remaining 40 images.  The training batch \(\K_i\) contains all four training images and is the same in every iteration \(i\).   Hence, we do not apply mini-batching in these experiments and the only changes between iterations are the neurons that are removed in the \texttt{ModifyNetwork} procedure in Algorithm 1.

We conduct experiments on two fully connected network architectures (Network 10 and 12 in Table~\ref{tab:networks}) having three and five neurons in the hidden layer, respectively.  In each iteration, we randomly remove five neurons (20~\%) from the input layer and $\ndrop=0, 1$, and $2$ neurons in the middle, hidden layer.  Output neurons are never removed.  The learning rate is varied using the values  $\learnrate=0.01, 0.02, 0.05, 0.1$, and $0.5$, while \(\beta=0.1\) is kept constant.   The number of iterations was kept constant for all experiments at \(M=10\).

For each combination of dropout hyperparameters, 100 independent training and testing runs were performed.  The settings for the number of replicas and Monte Carlo steps in each annealing calculation were identical to those described in Section~\ref{sec:setup}.  During the intermediate dropout iterations, the temperature parameters were fixed as \(T^{-1}_{\rm max}=0.2\) and \(T^{-1}_{\rm min} = 8.6\), but for the final training phase, the temperature parameters were optimised using the Nelder-Mead method~\cite{nelder1965}.  For comparison, 100 additional tests were conducted for each network without applying dropout (i.e., without iterative training) to serve as the baseline.
The results of the experiments are summarised in Table~\ref{tab:dropout-results-summary}, which lists for each set of experiments the minimum, maximum, mean, and median test accuracy as well as the mean training accuracy.

\begin{table}[htbp]
\caption{\label{tab:dropout-results-summary}Summary statistics for test and training accuracy for varying the learning rate \(\learnrate\) and the number of neurons dropped in the hidden layer \(n_{drop}\) for Network~10 and Network~12. Each entry represents results obtained from 100 independent training and inference runs on the four-class image classification task.  The table lists the minimum, maximum, mean, and median test accuracies, the corresponding mean training accuracies, and the mean fraction of unsatisfied constraints \eqref{eq:training-linear-substitution} and \eqref{eq:training-activation-constraint2}.}
\centering
\resizebox{\linewidth}{!}{
\begin{tabular}{ccccccccc}
\toprule
Network & $\learnrate$ & $\ndrop$ & \multicolumn{4}{c}{Test accuracy} & Training accuracy & Unsatisfied constraints\\[0pt]
\cline{4-7}
        &              &          & min. & max. & mean & median       & mean     & mean fraction (\%)\\[0pt]
\midrule
10 & 0.00 & - & 0.400 & 0.825 & 0.574 & 0.575 & 1.000 & 0.00\\[0pt]
10 & 0.01 & 0 & 0.375 & 0.850 & 0.583 & 0.575 & 1.000 & 0.00\\[0pt]
10 & 0.01 & 1 & 0.350 & 0.800 & 0.587 & 0.600 & 1.000 & 0.00\\[0pt]
10 & 0.01 & 2 & 0.325 & 0.875 & 0.577 & 0.575 & 1.000 & 0.00\\[0pt]
10 & 0.02 & 0 & 0.375 & 0.900 & 0.592 & 0.575 & 0.993 & 0.07\\[0pt]
10 & 0.02 & 1 & 0.250 & 0.825 & 0.568 & 0.575 & 0.983 & 0.16\\[0pt]
10 & 0.02 & 2 & 0.325 & 0.850 & 0.586 & 0.575 & 1.000 & 0.00\\[0pt]
10 & 0.05 & 0 & 0.350 & 0.825 & 0.575 & 0.550 & 0.960 & 0.55\\[0pt]
10 & 0.05 & 1 & 0.300 & 0.825 & 0.546 & 0.537 & 0.815 & 1.84\\[0pt]
10 & 0.05 & 2 & 0.300 & 0.800 & 0.578 & 0.550 & 0.965 & 0.32\\[0pt]
10 & 0.10 & 0 & 0.325 & 0.875 & 0.590 & 0.575 & 0.922 & 0.91\\[0pt]
10 & 0.10 & 1 & 0.275 & 0.800 & 0.540 & 0.550 & 0.777 & 2.48\\[0pt]
10 & 0.10 & 2 & 0.275 & 0.850 & 0.576 & 0.575 & 0.907 & 1.07\\[0pt]
10 & 0.50 & 0 & 0.225 & 0.875 & 0.571 & 0.575 & 0.780 & 2.20\\[0pt]
10 & 0.50 & 1 & 0.200 & 0.750 & 0.493 & 0.500 & 0.670 & 3.52\\[0pt]
10 & 0.50 & 2 & 0.275 & 0.825 & 0.537 & 0.525 & 0.723 & 2.77\\[0pt]
\midrule
12 & 0.00 & - & 0.275 & 0.875 & 0.542 & 0.550 & 1.000 & 0.00\\[0pt]
12 & 0.01 & 0 & 0.325 & 0.800 & 0.546 & 0.550 & 1.000 & 0.00\\[0pt]
12 & 0.01 & 1 & 0.400 & 0.775 & 0.579 & 0.600 & 1.000 & 0.00\\[0pt]
12 & 0.01 & 2 & 0.325 & 0.875 & 0.545 & 0.537 & 1.000 & 0.00\\[0pt]
12 & 0.02 & 0 & 0.275 & 0.775 & 0.551 & 0.550 & 0.960 & 0.28\\[0pt]
12 & 0.02 & 1 & 0.250 & 0.775 & 0.529 & 0.525 & 0.963 & 0.25\\[0pt]
12 & 0.02 & 2 & 0.325 & 0.875 & 0.573 & 0.575 & 0.995 & 0.06\\[0pt]
12 & 0.05 & 0 & 0.300 & 0.775 & 0.524 & 0.525 & 0.895 & 0.96\\[0pt]
12 & 0.05 & 1 & 0.250 & 0.825 & 0.518 & 0.525 & 0.848 & 1.25\\[0pt]
12 & 0.05 & 2 & 0.250 & 0.800 & 0.534 & 0.550 & 0.897 & 0.96\\[0pt]
12 & 0.10 & 0 & 0.250 & 0.775 & 0.492 & 0.500 & 0.765 & 2.24\\[0pt]
12 & 0.10 & 1 & 0.300 & 0.775 & 0.499 & 0.500 & 0.790 & 2.18\\[0pt]
12 & 0.10 & 2 & 0.300 & 0.775 & 0.551 & 0.575 & 0.873 & 1.25\\[0pt]
12 & 0.50 & 0 & 0.175 & 0.750 & 0.496 & 0.500 & 0.682 & 5.35\\[0pt]
12 & 0.50 & 1 & 0.250 & 0.725 & 0.463 & 0.450 & 0.568 & 5.93\\[0pt]
12 & 0.50 & 2 & 0.350 & 0.875 & 0.589 & 0.600 & 0.810 & 5.00\\[0pt]
\bottomrule
\end{tabular}
}
\end{table}

\begin{figure}
  \centering

  \begin{subfigure}{0.48\textwidth}
    \centering
    \includegraphics[width=\linewidth]{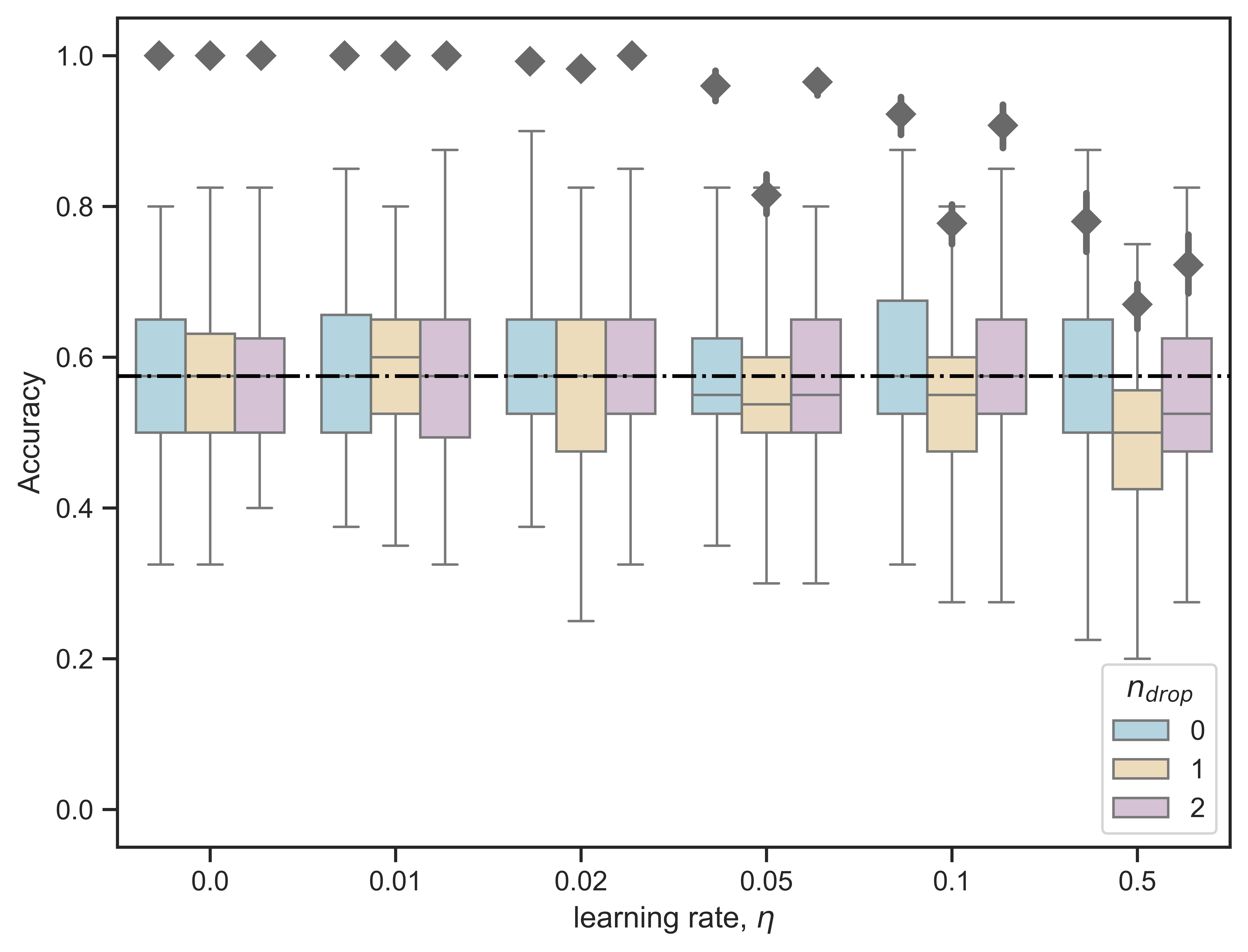}
    \caption{Network 10.}
    \label{fig:dropout-network10}
  \end{subfigure}
  %% \hspace{1cm}
  \hfill
  \begin{subfigure}{0.48\textwidth}
    \centering
    \includegraphics[width=\linewidth]{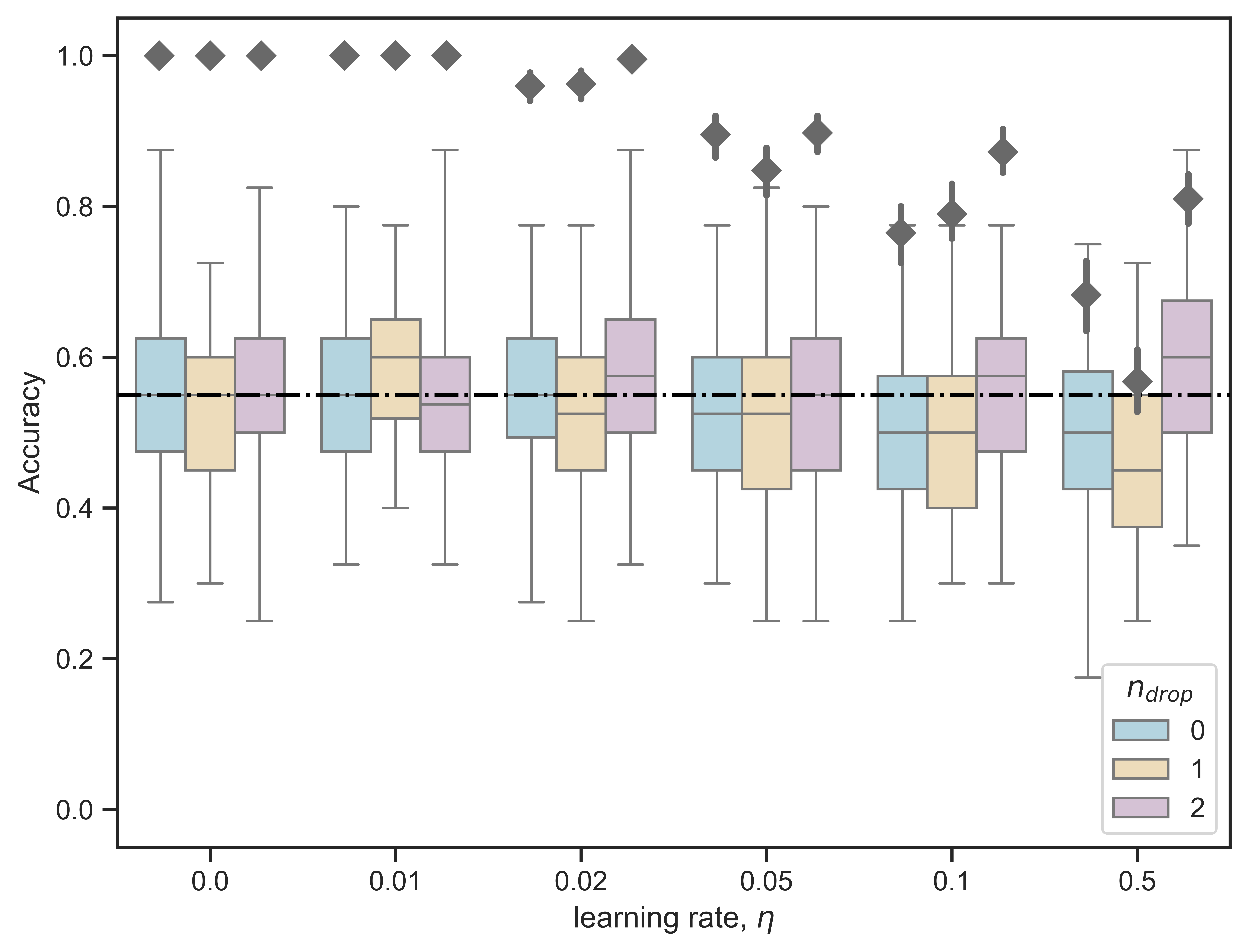}
    \caption{Network 12.}
    \label{fig:dropout-network12}
  \end{subfigure}

\caption{Test accuracy distributions for various learning rates $\eta$ and dropout levels $n_{drop}$ for Network 10~\subrefp{fig:dropout-network10} and Network 12~\subrefp{fig:dropout-network12}, for 100 independent training and inference runs.  Each box shows the first quartile, median, and third quartile, while the whiskers indicate the minimum and maximum.  The horizontal dashed line indicates the median test accuracy obtained without applying dropout.  Grey squares and vertical lines represent the mean and range of training accuracy, respectively.}
\label{fig:dropout-results}
\end{figure}

Figure~\ref{fig:dropout-results} plots the test accuracy by learning rate and dropout rate for each of the two networks considered.    The effect of the dropout algorithm on both networks is rather modest in this experiment.
For Network~12, slight improvements are occasionally observed when both \(\learnrate\) and \(\ndrop\) increase.  In particular, with a learning rate of \(\learnrate = 0.5\) and \(\ndrop = 2\), the mean test accuracy improves from 54\% (without dropout) to 59\%.  In contrast, Network 10 exhibits a slight decrease in accuracy under the same conditions.
This degradation is likely due to the small number of neurons (three) in the hidden layer, making the architecture less suitable for dropout.  In addition, the training accuracy generally tends to decrease as \(\learnrate\) increases.
This behaviour might be attributed to the additional penalty term \(H_{ext}\) becoming dominant over \(H_1\) and \(H_2\) in the final annealing run, which results in infeasible solutions for constraints~\eqref{eq:training-activation-constraint2} and \eqref{eq:training-linear-substitution}.
However, the reduction in test accuracy with increasing \(\learnrate\) remains relatively small.
This may serve as evidence that the dropout process improves the generalisation ability of the network during training.
To further verify the effects of model size reduction by batching and the generalisation benefits of the dropout, additional experiments using larger datasets and network architectures will be required.

\section{Conclusion}
\label{sec:conclusion}

QUBO-based formulations have recently enabled the training of BNNs on Ising machines.  In this work we show that the formulation can be extended to an arbitrary network topology, which is not necessarily layered or acyclic.  We conduct computational experiments on a small four-class image classification problem on a digital Ising machine and introduce two novel methods for regularising the training problem.
One is a quadratic penalty that maximises the sum of neuron margins.  Computational results show that this method significantly improves test accuracy and produces a clear positive correlation between margin size and generalisation.  The other is a dropout-like algorithm that randomly modifies the network topology during training.  The results do not show performance improvements in our experiments, but test accuracy stays robust even at higher dropout rates, indicating a modest generalisation benefit.

Despite these results, several limitations remain.
The training approach is restricted by the algorithm's ability to find solutions that satisfy all activation constraints.  As the size of the network and the training data set increases, this becomes increasingly difficult and impacts performance.  Hence, the experimental results presented in this paper are limited to a relatively small, constructed data set.  Further studies are necessary to confirm the methods' applicability to larger and more realistic data sets along with potential algorithmic and model improvements.

Future research should seek to improve performance, scalability, and applicability through algorithmic and model enhancements.
Firstly, when formulating the training problem as a QUBO problem, an energy landscape that can be efficiently explored by the annealing algorithm is desired for improving performance~\cite{dobrynin2024}. From this perspective, it is important to determine appropriate scaling factors $\alpha$ and $\gamma$, and to refine the formulation itself based on both theoretical and experimental analysis.
Secondly, further improvements and customisations of the annealing algorithm, such as enhancing the treatment of constraints~\cite{kameyama2024}, may aid in finding high-quality solutions to the training problem.
Thirdly, the assumptions in this study that model parameters are binary and that labelled training data are always available may be overly restrictive in practice, potentially limiting performance and applicability. However, multi-bit weights can be introduced with only minor modifications by allowing parallel bit-level connections, while extending the formulation to unsupervised learning on unlabelled data suggests a promising direction for future investigation.
Finally, exploring tasks that may benefit from cyclic network structures enabled by the model presented herein constitute a further potential area of research.

\section*{Declaration of generative AI and AI-assisted technologies in the manuscript preparation process}
During the preparation of this work the authors used ChatGPT in order to improve language and readability. The authors have reviewed and edited the provided suggestions as needed and take full responsibility for the content of the published article.

\section*{Acknowledgements}
We thank Alireza Ahrabian, Takuya Okuyama, Kento Hasegawa, Kazuo Ono, and Masanao Yamaoka for fruitful discussions and suggestions and for providing the computational environment for our experiments.

\bibliographystyle{elsarticle-num} 
\bibliography{references.bib}

\end{document}